\documentclass[11pt]{amsart}

\usepackage{graphics}
\usepackage{color}
\usepackage[a4paper, margin=3cm]{geometry}

\usepackage{amssymb,enumerate}
\usepackage{amsmath}
\usepackage{bbm}           
\usepackage{bm}
\usepackage{eso-pic,graphicx}
\usepackage{tikz}
\usepackage{cite}
\usepackage{esint}
\usepackage[colorlinks=true, pdfstartview=FitV, linkcolor=blue, citecolor=blue, urlcolor=blue]{hyperref}

\usepackage{graphicx}
\usepackage{pslatex}
\usepackage{amsmath,amsfonts}
\usepackage{rotating}
\usepackage{amssymb}
\usepackage{verbatim}
\usepackage{rotating}
\usepackage{mathrsfs}

\makeatletter
\def\Ddots{\mathinner{\mkern1mu\raise\p@
\vbox{\kern7\p@\hbox{.}}\mkern2mu
\raise4\p@\hbox{.}\mkern2mu\raise7\p@\hbox{.}\mkern1mu}}
\makeatother

\def\XXint#1#2#3{{\setbox0=\hbox{$#1{#2#3}{\int}$}
\vcenter{\hbox{$#2#3$}}\kern-.5\wd0}}

\begin{document}
\newtheorem{theorem}{Theorem}
\newtheorem{proposition}[theorem]{Proposition}
\newtheorem{conjecture}[theorem]{Conjecture}
\def\theconjecture{\unskip}
\newtheorem{corollary}[theorem]{Corollary}
\newtheorem{lemma}[theorem]{Lemma}
\newtheorem{claim}[theorem]{Claim}
\newtheorem{sublemma}[theorem]{Sublemma}
\newtheorem{observation}[theorem]{Observation}
\theoremstyle{definition}
\newtheorem{definition}{Definition}
\newtheorem{notation}[definition]{Notation}
\newtheorem{remark}[definition]{Remark}
\newtheorem{question}[definition]{Question}
\newtheorem{questions}[definition]{Questions}
\newtheorem{example}[definition]{Example}
\newtheorem{problem}[definition]{Problem}
\newtheorem{exercise}[definition]{Exercise}
 \newtheorem{thm}{Theorem}
 \newtheorem{cor}[thm]{Corollary}
 \newtheorem{lem}{Lemma}[section]
 \newtheorem{prop}[thm]{Proposition}
 \theoremstyle{definition}
 \newtheorem{dfn}[thm]{Definition}
 \theoremstyle{remark}
 \newtheorem{rem}{Remark}
 \newtheorem{ex}{Example}
 \numberwithin{equation}{section}
\def\C{\mathbb{C}}
\def\R{\mathbb{R}}
\def\Rn{{\mathbb{R}^n}}
\def\Rns{{\mathbb{R}^{n+1}}}
\def\Sn{{{S}^{n-1}}}
\def\M{\mathbb{M}}
\def\N{\mathbb{N}}
\def\Q{{\mathbb{Q}}}
\def\Z{\mathbb{Z}}
\def\F{\mathcal{F}}
\def\L{\mathcal{L}}
\def\S{\mathcal{S}}
\def\supp{\operatorname{supp}}
\def\essi{\operatornamewithlimits{ess\,inf}}
\def\esss{\operatornamewithlimits{ess\,sup}}

\numberwithin{equation}{section}
\numberwithin{thm}{section}
\numberwithin{theorem}{section}
\numberwithin{definition}{section}
\numberwithin{equation}{section}

\def\earrow{{\mathbf e}}
\def\rarrow{{\mathbf r}}
\def\uarrow{{\mathbf u}}
\def\varrow{{\mathbf V}}
\def\tpar{T_{\rm par}}
\def\apar{A_{\rm par}}

\def\reals{{\mathbb R}}
\def\torus{{\mathbb T}}
\def\scriptm{{\mathcal T}}
\def\heis{{\mathbb H}}
\def\integers{{\mathbb Z}}
\def\z{{\mathbb Z}}
\def\naturals{{\mathbb N}}
\def\complex{{\mathbb C}\/}
\def\distance{\operatorname{distance}\,}
\def\support{\operatorname{support}\,}
\def\dist{\operatorname{dist}\,}
\def\Span{\operatorname{span}\,}
\def\degree{\operatorname{degree}\,}
\def\kernel{\operatorname{kernel}\,}
\def\dim{\operatorname{dim}\,}
\def\codim{\operatorname{codim}}
\def\trace{\operatorname{trace\,}}
\def\Span{\operatorname{span}\,}
\def\dimension{\operatorname{dimension}\,}
\def\codimension{\operatorname{codimension}\,}
\def\nullspace{\scriptk}
\def\kernel{\operatorname{Ker}}
\def\ZZ{ {\mathbb Z} }
\def\p{\partial}
\def\rp{{ ^{-1} }}
\def\Re{\operatorname{Re\,} }
\def\Im{\operatorname{Im\,} }
\def\ov{\overline}
\def\eps{\varepsilon}
\def\lt{L^2}
\def\diver{\operatorname{div}}
\def\curl{\operatorname{curl}}
\def\etta{\eta}
\newcommand{\norm}[1]{ \|  #1 \|}
\def\expect{\mathbb E}
\def\bull{$\bullet$\ }

\def\blue{\color{blue}}
\def\red{\color{red}}

\def\xone{x_1}
\def\xtwo{x_2}
\def\xq{x_2+x_1^2}
\newcommand{\abr}[1]{ \langle  #1 \rangle}

\newcommand{\Norm}[1]{ \left\|  #1 \right\| }
\newcommand{\set}[1]{ \left\{ #1 \right\} }
\newcommand{\ifou}{\raisebox{-1ex}{$\check{}$}}
\def\one{\mathbf 1}
\def\whole{\mathbf V}
\newcommand{\modulo}[2]{[#1]_{#2}}
\def \essinf{\mathop{\rm essinf}}
\def\scriptf{{\mathcal F}}
\def\scriptg{{\mathcal G}}
\def\scriptm{{\mathcal M}}
\def\scriptb{{\mathcal B}}
\def\scriptc{{\mathcal C}}
\def\scriptt{{\mathcal T}}
\def\scripti{{\mathcal I}}
\def\scripte{{\mathcal E}}
\def\scriptv{{\mathcal V}}
\def\scriptw{{\mathcal W}}
\def\scriptu{{\mathcal U}}
\def\scriptS{{\mathcal S}}
\def\scripta{{\mathcal A}}
\def\scriptr{{\mathcal R}}
\def\scripto{{\mathcal O}}
\def\scripth{{\mathcal H}}
\def\scriptd{{\mathcal D}}
\def\scriptl{{\mathcal L}}
\def\scriptn{{\mathcal N}}
\def\scriptp{{\mathcal P}}
\def\scriptk{{\mathcal K}}
\def\frakv{{\mathfrak V}}
\def\C{\mathbb{C}}
\def\D{\mathcal{D}}
\def\R{\mathbb{R}}
\def\Rn{{\mathbb{R}^n}}
\def\rn{{\mathbb{R}^n}}
\def\Rm{{\mathbb{R}^{2n}}}
\def\r2n{{\mathbb{R}^{2n}}}
\def\Sn{{{S}^{n-1}}}
\def\M{\mathbb{M}}
\def\N{\mathbb{N}}
\def\Q{{\mathcal{Q}}}
\def\Z{\mathbb{Z}}
\def\F{\mathcal{F}}
\def\L{\mathcal{L}}
\def\G{\mathscr{G}}
\def\ch{\operatorname{ch}}
\def\supp{\operatorname{supp}}
\def\dist{\operatorname{dist}}
\def\essi{\operatornamewithlimits{ess\,inf}}
\def\esss{\operatornamewithlimits{ess\,sup}}
\def\dis{\displaystyle}
\def\dsum{\displaystyle\sum}
\def\dint{\displaystyle\int}
\def\dfrac{\displaystyle\frac}
\def\dsup{\displaystyle\sup}
\def\dlim{\displaystyle\lim}
\def\bom{\Omega}
\def\om{\omega}
\author[F. Liu]{Feng Liu}
\address{Feng Liu:
        College of Mathematics and System Science\\
        Shandong University of Science and Technology\\
        Qingdao, Shandong 266590\\
        People's Republic of China}
\email{FLiu@sdust.edu.cn}

\author[Y. Lian]{Yongzhi Lian}
\address{Yongzhi Lian:
        College of Mathematics and System Science\\
        Shandong University of Science and Technology\\
        Qingdao, Shandong 266590\\
        People's Republic of China}
\email{18764016097@163.com}

\author{Yanzhen Wang$^{*}$}
\address{Yanzhen Wang:
        College of Mathematics and System Science\\
        Shandong University of Science and Technology\\
        Qingdao, Shandong 266590\\
        People's Republic of China}
\email{wangyanzhen825@163.com}


\thanks{$^{*}$ Corresponding author, e-mail address: wangyanzhen825@163.com}

\date{\today}
\title[maximal operators on finite connected graphs]
{\bf Sharp $L^p$ and variation inequalities for maximal operators on finite connected graphs}

\begin{abstract} Let $G=(V_G,E_G)$ be a finite connected graph (here $V_G$
and $E_G$ denote the set of vertices and edges of $G$ respectively). The 
primary aim of this paper is to study the sharp variation and $L^p$ 
inequalities of the centered Hardy--Littlewood maximal operator $M_G$ 
defined on $G$ and its fractional variant $M_{\alpha,G}\,(0<\alpha<1)$). 
Set
$${\bf C}_{\alpha,G}^{p,q}=\sup\limits_{{\rm Var}_{p,G}(f)>0}\frac{{\rm Var}_{q,G}(M_{\alpha,G}f)}{{\rm Var}_{p,G}(f)},$$
where ${\rm Var}_{p,G}(f)$ denotes the $p$-variation of $f$ on $G$. We
denote by $K_n$ (resp., $S^n$) the complete (resp., star) graph with
$n\geq4$ vertices. When $\alpha=0$, we prove ${\bf C}_{\alpha,S_n}^{p,p}
=1-\frac{1}{n}$ for $0<p<1$. This solve a question originally posed by 
Liu and Xue in 2021. In the fractional case $0<\alpha<1$, we prove 
${\bf C}_{\alpha,S_n}^{1,q}=1-n^{-\alpha}$ and ${\bf C}_{\alpha,K_n}^{1,q}
=(1-n^{\alpha-1})(n-1)^{-\alpha}$ for $q=\frac{1}{1-\alpha}$, which provide
some examples focusing on sharp variation inequalities of fractional maximal
operators. The sharp variation inequalities of uncentered fractional maximal
operators are also established. In addition, the sharp $L^p$ inequalities of
uncentered fractional maximal operator $\widetilde{M}_{\alpha,G}$ are also
studied. More precisely, we obtain lower bounds and upper bounds for $C_{q,G}$
such that the inequality $\|\widetilde{M}_{\alpha,G}f\|_q\leq C_{q,G}\|f\|_p$
holds and the complete graphs $K_n$ and star graphs $S_n$ are two extreme
graphs that reach the lower bound and upper bound respectively.
\end{abstract}\maketitle

\section{Introduction}\label{S1}

It was well known that maximal operator is one of the most important
operators in harmonic analysis, which has a wide variety of applications,
such as harmonic analysis, ergodic theory, approximation theory, PDEs and
so on. The study of sharp inequalities of various maximal operators on
Lebesgue spaces and other function spaces has been always a challenging
topic of current research, one may consult \cite{Ald,GM3,GMM,Me1,Me2}
for the sharp inequalities of Hardy--Littlewood maximal operators on
Lebesgue spaces, \cite{AL,CS,Ra} for the sharp variation inequalities
of various maximal operators. In this paper we focus on the sharp $L^p$
inequalities and variation inequalities for Hardy--Littlewood maximal
operators on finite graphs and their fractional variants, both in the
centered and uncentered versions.

\subsection{Sharp inequalities of maximal operators on graphs on Lebesgue spaces}
Throughout this paper, we denote $\mathbb{N}=\{0,1,\ldots\}$ and
$\mathbb{R}^{+}=[0,\infty)$. Let $G=(V_G,E_G)$ be an undirected
combinatorial graph with the set of vertices $V_G$ and the set of
edges $E_G$. Two vertices $x,\,y\in V_G$ are called neighbors if
they are connected by an edge $x\sim y\in E_G$. For any $v\in V_G$
let $N_G(v)$ be the set of neighbors of $v$. The graph $G$ is finite
if $|V_G|<\infty$ for any $v\in V_G$. We say that $G$ is called
connected if for any distinct $x,\,y\in V_G$, there exists a finite
sequence of vertices $\{x_i\}_{i=0}^k$, $k\in\mathbb{N}\setminus\{0\}$,
such that $x=x_0\sim x_1\sim\cdots\sim x_k=y$. We say that $G$ is
infinite if $|V_G|=+\infty$. Two graphs $G_1=(V_{G_1},E_{G_1})$ and
$G_2=(V_{G_2},E_{G_2})$ are isomorphic if there is a bijection of
the vertices $\theta:V_{G_1}\rightarrow V_{G_2}$ such that
$v\sim w\in E_{G_1}$ if and only if $\theta(v)\sim\theta(w)\in E_{G_2}$
for any vertices $v,\,w\in V_{G_1}$.

In what follows, we assume that $G=(V_G,E_G)$ is a finite connected
graph with the set of vertices $V_G$ and the set of edges $E_G$. Let
$d_G$ be the metric induced by the edges in $E_G$, i.e. for given
$u,v\in V_G$, the distance $d_G(u,v)$ is the number of edges in a
shortest path connecting $u$ and $v$. Let $B_G(v,r)$ be the ball
centered at $v$, with radius $r$ on the graph $G$, i.e.
$$B_G(v,r)=\{u\in V_G:d_G(u,v)\leq r\}.$$
For example,
\[
  B_G(v,r)= \left\{
\begin{array}{ll}
\{v\},&\ {\rm if}\ 0\leq r<1;\\
\{v\}\cup N_G(v),&\ {\rm if}\ 1\leq r<2.
\end{array}
\right.
\]
For any subset $A\subset V_G$, the notation $|A|$ represents the cardinality
of $A$. In order to establish our main results, let us introduce two particular
graphs.

$\bullet$ (Complete graph $K_n$). Let $K_n=(V_{K_n},E_{K_n})$ be the
complete graph with $n\geq2$ vertices, i.e. $N_{K_n}(v)=V_{K_n}
\setminus\{v\}$ for all $v\in V_{K_n}$.

$\bullet$ (Star graph $S_n$). Let $S_n=(V_{S_n},E_{S_n})$ be the star
graph with $n\geq2$ vertices, i.e. there exists an unique $v\in V_{S_n}$
such that $N_{S_n}(v)=V_{S_n}\setminus\{v\}$ and $N_{S_n}(w)=\{v\}$ for
all $w\in V_{S_n}\setminus\{v\}$.

We now introduce the maximal operators on graphs.

\begin{definition}\label{def1.1}
{\bf (Maximal operators on graphs)} Let $G=(V_G,E_G)$ and $0\leq\alpha<1$.
For a function $f:V_G\rightarrow\mathbb{R}$, the centered fractional
maximal operator on $G$ is defined by
$$M_{\alpha,G}f(v)=\sup\limits_{r\in\mathbb{N}}\frac{1}{|B_G(v,r)|^{1-\alpha}}\sum\limits_{w\in B_G(v,r)}|f(w)|,\ \ \ v\in V_G.$$
The uncentered version is given by
$$\widetilde{M}_{\alpha,G}f(v)=\sup\limits_{v\in B_G(w,r),\,r\in\mathbb{N}}\frac{1}{|B_G(w,r)|^{1-\alpha}}\sum\limits_{u\in B_G(w,r)}|f(u)|,\ \ \ v\in V_G.$$
When $\alpha=0$, the operator $M_{\alpha,G}$ (resp., $\widetilde{M}_{\alpha,G}$)
reduces to the centered (resp., uncentered) Hardy--Littlewood maximal operator
on $G$, which is denoted by $M_G$ (resp., $\widetilde{M}_G$).
\end{definition}

The Hardy--Littlewood maximal operator on graphs has been studied by many
authors (see \cite{BM,CMS,GM,GM2,KP,ST,ST2}). This topic originated from
Kor\'{a}nyi and Picardello's work \cite{KP} who used the Hardy--Littlewood
maximal operators on graphs to study the boundary behaviour of
eigenfunctions of the Laplace operator on trees. Later on, Cowling et al.
\cite{CMS} studied the Hardy--Littlewood maximal operators on homogeneous
trees. Badr and Martell \cite{BM} established the weighted norm inequalities 
for the Hardy--Littlewood maximal operators on infinite graphs. Soria and 
Tradacete \cite{ST2} further investigated some different geometric properties 
on infinite graphs, related to the weak-type boundedness of the 
Hardy--Littlewood maximal operator.

In the finite graph setting, Soria and Tradacete \cite{ST} firstly studied
the $L^p$-norm for the Hardy--Littlewood maximal operators on finite
complete graph and star graph. Recently, Xiang and Zhao \cite{XZ} studied
the operator norm of centered fractional maximal operator on Lebesgue
spaces.

We now introduce one definition.

\begin{definition}\label{def1.2}
{\bf (Lebesgue space on graphs)} Let $G=(V_G,E_G)$ and $0<p\leq\infty$. The
$L^p(V_G)$ space is the set of all functions $f:V_G\rightarrow\mathbb{R}$
satisfying $\|f\|_p<\infty$, where
\[
  \|f\|_p= \left\{
\begin{array}{ll}
\displaystyle\Big(\sum\limits_{v\in V_G}|f(v)|^p\Big)^{1/p},&\ {\rm if}\ 0<p<\infty;\\
\sup\limits_{v\in V_G}|f(v)|,&\ {\rm if}\ p=\infty.
\end{array}
\right.
\]
\end{definition}
Let $0\leq\alpha<1$ and $0\leq p\leq q\leq\infty$. We define
$$\|M_{\alpha,G}\|_{L^p\rightarrow L^q}=\sup\limits_{f:V_G\rightarrow\mathbb{R},\,f\neq0}\frac{\|M_{\alpha,G}f\|_q}{\|f\|_p},\ \ \
\|\widetilde{M}_{\alpha,G}\|_{L^p\rightarrow L^q}=\sup\limits_{f:V_G\rightarrow\mathbb{R},\,f\neq0}\frac{\|\widetilde{M}_{\alpha,G}f\|_q}{\|f\|_p}.$$

The main results of \cite{ST,XZ} can be formulated as follows.

\medskip

\quad\hspace{-20pt}{\bf Theorem A} {\rm (\cite{ST,XZ})} {\it Let
$G=(V_G,E_G)$ be a finite connected graph with $n\geq4$ vertices,
and let $0\leq\alpha<1$ and $0<p\leq q\leq1$, then
$$\Big(1+\frac{n-1}{n^{q(1-\alpha)}}\Big)^{{1}/{q}}\leq\|M_{\alpha,G}\|_{L^p\rightarrow L^q}\leq\Big(1+\frac{n-1}{2^{q(1-\alpha)}}\Big)^{{1}/{q}}.$$
Moreover,
\begin{enumerate}[{\rm (i)}]
\item $\|M_{\alpha,G}\|_{L^p\rightarrow L^q}=\big(1+\frac{n-1}{n^{q(1-\alpha)}}\big)^{{1}/{q}}$
if and only if $G=K_n$;
\item $\|M_{\alpha,G}\|_{ L^p\rightarrow L^q}=\big(1+\frac{n-1}{2^{q(1-\alpha)}}\big)^{{1}/{q}}$
if and only if $G$ is isomorphic to $S_n$ ($n\geq3$).
\end{enumerate}}

It is worth noting that the operator norms of centered and uncentered
maximal operators on Lebesgue spaces are distinct in the European space.
However, the above claim doesn't seem necessary true. For example,
when $G=K_n$, we see that $M_{\alpha,K_n}=\widetilde{M}_{\alpha,K_n}$.
Going back to the star graph $S_n$, it is easy to see that
$M_{\alpha,S_n}f(v)\leq \widetilde{M}_{\alpha,S_n}f(v)$ for all
$v\in V_{S_n}$ and any function $f:V_{S_n}\rightarrow\mathbb{R}$. But
there does not exist a constant $C>0$ such that $\widetilde{M}_{\alpha,
S_n}f(v)\leq CM_{\alpha,S_n}f(v)$ for all $v\in V_{S_n}$ and any
function $f:V_{S_n}\rightarrow\mathbb{R}$. It makes the study of
$\|\widetilde{M}_{\alpha,G}\|_{L^p\rightarrow L^q}$ interesting.

Based on Theorem A, it is natural to ask the following question.

\begin{question}\label{que1.1}
Does the conclusion of Theorem A (ii) hold for $\widetilde{M}_{\alpha,S_n}$?
\end{question}

In this paper we shall prove the following result.

\begin{theorem}\label{thm1.1}
Let $G=(V_G,E_G)$ be a finite connected graph with $n\geq4$ vertices,
and let $0\leq\alpha<1$ and $0<p\leq q\leq1$, then
$$\Big(1+\frac{n-1}{n^{q(1-\alpha)}}\Big)^{{1}/{q}}\leq\|\widetilde{M}_{\alpha,G}\|_{L^p\rightarrow L^q}\leq\Big(1+\frac{n-1}{2^{q(1-\alpha)}}\Big)^{{1}/{q}}.\eqno(1.1)$$
Moreover,
\begin{enumerate}[{\rm (i)}]
\item $\|\widetilde{M}_{\alpha,G}\|_{L^p\rightarrow L^q}=\big(1+\frac{n-1}{n^{q(1-\alpha)}}\big)^{{1}/{q}}$
if and only if $G=K_n$;
\item $\|\widetilde{M}_{\alpha,G}\|_{L^p\rightarrow L^q}=\big(1+\frac{n-1}{2^{q(1-\alpha)}}\big)^{{1}/{q}}$
if and only if $G$ is isomorphic to $S_n$ ($n\geq3$).
\end{enumerate}
\end{theorem}

In addition, we shall establish the extremizers for $\|\widetilde{M}_{\alpha,
S_n}\|_{L^p\rightarrow L^p}$ with $1<p<\infty$.

\begin{theorem}\label{thm1.2}
Let $n\geq3$ and $S_n=(V_{S_n},E_{S_n})$ be the star graph centered at
the vertex $v_0\in V$. Let $0\leq\alpha<1$ and $1<p<\infty$. There
exists $f:V_{S_n}\rightarrow\mathbb{R}^+$ such that $f(v_0)=
\max_{j\in V_{S_n}}f(j)$, $f|_{V_{S_n}\setminus\{v_0\}}$ takes at
most three values, and
$$\|\widetilde{M}_{\alpha,S_n}\|_{L^p\rightarrow L^p}=\frac{\|\widetilde{M}_{\alpha,S_n} f\|_p}{\|f\|_p}.\eqno(1.2)$$
\end{theorem}
\begin{remark}\label{rem1.1}
\begin{enumerate}[{\rm (i)}]
\item In view of Theorem A and Theorem \ref{thm1.1}, we see that
the $L^p$-norms of uncentered Hardy--Littlewood maximal operator
on graphs and its fractional variant have the same lower bounds and
upper bounds as the centered cases. Moreover,
$$\|\widetilde{M}_{\alpha,S_n}\|_{L^p\rightarrow L^q}=\|M_{\alpha,S_n}\|_{L^p\rightarrow L^q},$$
which are very interesting.
\item In \cite{XZ}, the authors studied the extremizers for $\|M_{\alpha,
K_n}\|_{L^p\rightarrow L^p}$ and $\|M_{\alpha,S_n}\|_{L^p\rightarrow L^p}$.
Particularly, the authors \cite{XZ} proved that under the same conditions
of Theorem \ref{thm1.2}, there exists a function $f:V_{S_n}\rightarrow
\mathbb{R}^+$ such that $f(v_0)=\max_{j\in V_{S_n}}f(j)$, $f|_{V_{S_n}
\setminus\{v_0\}}$ takes at most two values, and
$$\|M_{\alpha,S_n}\|_{L^p\rightarrow L^p}=\frac{\|M_{\alpha,S_n} f\|_p}{\|f\|_p}.$$
Compare the centered case, the extremizer function in the uncentered case
is more relaxed.
\end{enumerate}
\end{remark}

\subsection{Sharp variation inequalities of maximal operators on graphs}
An active topic of current research in harmonic analysis is the investigation
on the variation inequalities of various maximal operators, which is closely
related to the Sobolev regularity theory of maximal operators. There are
lots of works concerning variation inequalities of various maximal operators.
For example, see \cite{AL,Ku} for the one-dimensional Hardy--Littlewood
maximal operator, \cite{Ra} for the one-dimensional non-tangential maximal
operator, \cite{CM} for the one-dimensional uncentered fractional maximal
operator, \cite{CS,CFS} for the heat flow maximal operator and the Poisson
maximal operator, \cite{Wei1} for the uncentered maximal characteristic
function, \cite{Wei2} for the dyadic maximal operator, \cite{Wei3} for the
uncentered maximal operator with respect to cubes, \cite{BCHP,CM} for the
discrete Hardy--Littlewood maximal operators and their fractional variants,
among others.

Recently, the variation inequalities of maximal operators on graphs has
been studied by many authors (\cite{GKM,GM,LX}). Let us recall one
definition.

\begin{definition} [\textbf{${\rm BV}_{p}(G)$ space, \cite{LX}}] Let
$G=(V_G,E_G)$ and $0<p\leq\infty$. The space of bounded $p$-variation
functions is given by
$${\rm BV}_{p}(G):=\{f:V_G\rightarrow\mathbb{R};\ \|f\|_{{\rm BV}_{p}(G)}={\rm Var}_{p,G}(f)<\infty\},$$
where ${\rm Var}_{p,G}(f)$ represents the $p$-variation of $f$ defined by
\[
 {\rm Var}_{p,G}(f)= \left\{
\begin{array}{ll}
\displaystyle\Big(\sum\limits_{u\sim v\in E_G}|f(u)-f(v)|^p\Big)^{1/p},&\ {\rm if}\ 0<p<\infty;\\
\sup\limits_{u\sim v\in E_G}|f(u)-f(v)|,&\ {\rm if}\ p=\infty.
\end{array}
\right.
\]
We denote ${\rm BV}_{p}(G)={\rm BV}(G)$ and ${\rm Var}_{p,G}={\rm Var}_{G}$ 
when $p=1$.
\end{definition}

Let $0\leq\alpha<1$, $0<p,\,q<\infty$ and $G=(V_G,E_G)$. We define the
variation norm of maximal operators as follows:
$${\bf C}_{\alpha,G}^{p,q}=\sup\limits_{f:V_G\rightarrow\mathbb{R},\atop {\rm Var}_{p,G}(f)>0}\frac{{\rm Var}_{q,G}(M_{\alpha,G}f)}{{\rm Var}_{p,G}(f)},\ \
\widetilde{{\bf C}}_{\alpha,G}^{p,q}=\sup\limits_{f:V_G\rightarrow\mathbb{R},\atop {\rm Var}_{p,G}(f)>0}\frac{{\rm Var}_{q,G}(\widetilde{M}_{\alpha,G}f)}{{\rm Var}_{p,G}(f)}.$$
We denote
$${\bf C}_{\alpha,G}^{p,q}={\bf C}_G^{p,q},\ \ \ \widetilde{{\bf C}}_{\alpha,G}^{p,q}=\widetilde{{\bf C}}_G^{p,q},\ \ {\rm when}\ \alpha=0.$$

The variation inequality of maximal operators on graphs was firstly
studied by Liu and Xue \cite{LX} who observe that ${\bf C}_{K_3}^{p,p}=
{\bf C}_{S_3}^{p,p}=\frac{2}{3}$ for $0<p\leq1$. Then also proved that
$1-\frac{1}{n}<{\bf C}_{\alpha,K_n}^{p,p}<1$ and
$1-\frac{1}{n}<{\bf C}_{\alpha,S_n}^{p,p}<1$ for any $n\geq4$ and
$0<p\leq\infty$. Particularly, the above authors conjectured that
$${\bf C}_{K_n}^{p,p}={\bf C}_{S_n}^{p,p}=1-\frac{1}{n},\ \ {\rm if} \ n\geq4,\ \ 0<p\leq\infty.$$
In the case $G=K_n$, Madrid and his co-authors \cite{GKM,GM} proved
that if $n\geq3$ and $0<p<\infty$, then ${\bf C}_{K_n}^{p,p}=
1-\frac{1}{n}$. In the case $G=S_n$, Gonz\'{a}lez-Riquelme and Madrid
\cite{GM} showed that

\medskip

\quad\hspace{-20pt}{\bf Theorem B} {\rm (\cite{GM})} {\it Let $n=4$,
$0<p\leq1$ or $n\geq5$ and $\frac{1}{2}\leq p\leq1$, then
${\bf C}_{S_n}^{p,p}=1-\frac{1}{n}$.}

\medskip

Based on Theorem B, it is natural to ask the following question

\begin{question}\label{que1.2}
Whether ${\bf C}_{S_n}^{p,p}=1-\frac{1}{n}$ holds for $p\in(0,\frac{1}{2})$.
What happens when we consider the uncentered case?
\end{question}

In this paper we shall provide an affirmative answer to the above question.

\begin{theorem}\label{thm1.3}
Let $n\geq4$ and $0<p\leq1$. Then we have
$${\bf C}_{S_n}^{p,p}=1-\frac{1}{n}.\eqno(1.3)$$
\end{theorem}

\begin{theorem}\label{thm1.4}
Let $n\geq4$ and $0<p\leq1$. Then we have
$$\widetilde{{\bf C}}_{S_n}^{p,p}=\Big(2^{-p}+(n-2)\Big(\frac{n-2}{2n}\Big)^p\Big)^{1/p}.\eqno(1.4)$$
\end{theorem}

\begin{remark}\label{rem1.2}
Theorem \ref{thm1.3} answers a question posed by Liu and Xue \cite{LX}
and Gonz\'{a}lez-Riquelme and Madrid \cite{GM}. It should be pointed
out that the sharp variation constants of maximal operators on graph
$G$ can not characterize the properties of $G$. For example, it was
shown in \cite{GKM} that ${\bf C}_{G}^{1,1}=\frac{3}{4}$ if $G$
is a connected graph with $4$ vertices.
\end{remark}

Compare with the Hardy--Littlewood maximal operator, the study of the
variation inequalities of fractional maximal operators is few. When
$G=(V_G,E_G)$ with $V_G=\mathbb{Z}$ and $E_G=\{i\sim i+1:i\in\mathbb{Z}\}$,
$\widetilde{M}_{\alpha,G}$ reduces to the usual discrete uncentered
fractional maximal operator. Carneiro and Madrid \cite{CM} proved that
if ${\rm Var}_{G}(f)<\infty$ and $\widetilde{M}_{\alpha,G}f\not\equiv0$,
then
$${\rm Var}_q(\widetilde{M}_{\alpha,G}f)\leq 4^{1-\alpha}{\rm Var}_{G}(f),$$
where $q=\frac{1}{1-\alpha}$. However, the sharp constant
$\widetilde{{\bf C}}_{\alpha,G}^{1,q}$ is unknown. The variation inequality
of $M_{\alpha,G}$ is also unclear.

Based on the above, it is natural to ask whether one can obtain the sharp
$\frac{1}{1-\alpha}$-variation inequalities of $M_{\alpha,G}$ and
$\widetilde{M}_{\alpha,G}$ on finite graph $G$. More precisely, the above
variation inequality of $\widetilde{M}_\alpha$ stimulates the following
question.

\begin{question}\label{que1.3}
Let $0<\alpha<1$ and $q=\frac{1}{1-\alpha}$. What are ${\bf C}_{\alpha,G}^{1,q}$
and $\widetilde{{\bf C}}_{\alpha,G}^{1,q}$ when $G=K_n$ or $S_n$?
\end{question}

The above question can be answered by the following results.

\begin{theorem}\label{thm1.5}
Let $n\geq4$, $0<\alpha<1$ and $q=\frac{1}{1-\alpha}$. Then we have
$${\bf C}_{\alpha,K_n}^{1,q}=\widetilde{{\bf C}}_{\alpha,K_n}^{1,q}=(1-n^{\alpha-1})(n-1)^{-\alpha}.\eqno(1.5)$$
\end{theorem}

\begin{theorem}\label{thm1.6}
Let $n\geq4$, $0<\alpha<1$ and $q=\frac{1}{1-\alpha}$. Then we have
$${\bf C}_{\alpha,S_n}^{1,q}=1-n^{\alpha-1}.\eqno(1.6)$$
\end{theorem}

\begin{theorem}\label{thm1.7}
Let $n\geq4$, $0<\alpha<1$ and $q=\frac{1}{1-\alpha}$. Then we have
$$\widetilde{{\bf C}}_{\alpha,S_n}^{1,q}=((1-2^{\alpha-1})^p+(n-2)(2^{\alpha-1}-n^{\alpha-1})^p)^{1/p}.\eqno(1.7)$$
\end{theorem}

\subsection{Outline of this paper and notation}We now describe the outline
of the remainder of this paper. In Section \ref{S2} we prove Theorems
\ref{thm1.1} and \ref{thm1.2}. These proofs are motivated by the proofs
of \cite[Theorems 2.1 and 3.2]{XZ}. In Section \ref{S3} we prove Theorems
\ref{thm1.3} and \ref{thm1.4}. The proofs of Theorems \ref{thm1.5}--
\ref{thm1.7} will be given in Section \ref{S4}. In order to prove our main
results, the following Kronecker delta function $\delta_k\,(k\geq1)$ are
very useful, which are defined by
\[
\delta_k(i) =
\begin{cases}
1, & i = k, \\
0, & i \neq k.
\end{cases}
\]

\bigskip

\section{Operator norms of $\widetilde{M}_{\alpha,G}$ from $L^p$ to $L^q$--Proofs of Theorems \ref{thm1.1} and \ref{thm1.2}}\label{S2}
This section is devoted to presenting the proofs of Theorems \ref{thm1.1}
and \ref{thm1.2}.

\begin{proof}[{\rm Proof of Theorem \ref{thm1.1}}]
Applying the arguments similar to the proofs of \cite[Theorem 2.1]{XZ},
we can get (1.1). Theorem \ref{thm1.1} (i) follows from Theorem A (i).
Next we prove part (ii). This proof will be divided into two steps:

{\bf Step 1.} Proof of
$\|\widetilde{M}_{\alpha,S_n}\|_{L^p\rightarrow L^q}=
(1+\frac{n-1}{2^{q(1-\alpha)}})^{{1}/{q}}$.

Without loss of generality we may assume $S_n=(V_{S_n},E_{S_n})$, where
$V_{S_n}=\{1,\ldots,n\}$ and $E_{S_n}=\{1\sim2,\ldots,1\sim n\}$. Clearly,
\[
\widetilde{M}_{\alpha,S_n}\delta_1(j) =
\begin{cases}
1, & j=1, \\
2^{\alpha-1}, & j\in\{2,\ldots,n\}.
\end{cases}
\]
For any $2\leq i\leq n$, we have
\[
\widetilde{M}_{\alpha,S_n}\delta_i(j) =
\begin{cases}
2^{\alpha-1}, & j=1, \\
1, & j=i, \\
n^{\alpha-1}, & j\in\{2,\ldots,n\}\setminus\{i\}.
\end{cases}
\]
Applying \cite[Lemma 2.1]{XZ}, one gets
$$\|\widetilde{M}_{\alpha,S_n}\|_{L^p\rightarrow L^q}=\max\limits_{j\in V_{S_n}}\|\widetilde{M}_{\alpha,S_n}\delta_j\|_{q}=\Big(1+\frac{n-1}{2^{q(1-\alpha)}}\Big)^{{1}/{q}}.$$

{\bf Step 2.} Prove of $G$ is isomorphic to $S_n$
when $n\geq3$ and $\|\widetilde{M}_{\alpha,G}\|_{L^p\rightarrow L^q}=
(1+\frac{n-1}{2^{q(1-\alpha)}})^{{1}/{q}}$.

Let $n\geq3$ and $\|\widetilde{M}_{\alpha,G}\|_{L^p\rightarrow L^q}=
(1+\frac{n-1}{2^{q(1-\alpha)}})^{{1}/{q}}$. Given $f:V_{S_n}\rightarrow
\mathbb{R}^+$ with $\|f\|_1 = 1$, we get by H\"{o}lder's inequality that
$\|f\|_p \leq n^{\frac{1}{p}-1}\|f\|_1$. By \cite[Lemma 2.1]{XZ}, there
exists $v_0\in V_{G}$ such that
$$\|\widetilde{M}_{\alpha,G}\delta_{i_0}\|_q^q=1+\frac{n-1}{2^{q(1-\alpha)}}.\eqno(2.1)$$
Since $\widetilde{M}_{\alpha,G}\delta_{v_0}(v_0)=1$ and $\widetilde{M}_{\alpha,G}
\delta_{v_0}(w)\leq2^{\alpha-1}$ for all $w\neq v_0$. In view of (2.1), one gets
$$\widetilde{M}_{\alpha,G}\delta_{v_0}(w)=2^{\alpha-1},\ \ w\neq v_0.\eqno(2.2)$$
For any $w\in V_G\setminus\{v_0\}$. By (2.2), there exists a ball
$B_G(u_w,r_w)=\{v_0,w\}$ such that
$$\widetilde{M}_{\alpha,G}\delta_{v_0}(w)=\frac{1}{|B_G(u_w,r_w)|^{1-\alpha}}\sum\limits_{x\in B_G(u_w,r_w)}\delta_{v_0}(x)=\frac{1}{|B_G(u_w,r_w)|^{1-\alpha}}.$$
Since $G$ is a finite connected graph and $B_G(u_w,r_w)=\{v_0,w\}$, then
$r_w=d_G(v_0,w)=1$ and $u_w=v_0$ or $u_w=w$. If there exists
$w_0\in V_G\setminus\{v_0\}$ such that $u_{w_0}=v_0$. Then $|N_{G}(v_0)|=\{w_0\}$.
Let $w_1\in V_G\setminus\{v_0\}$ with $w_1\neq w_0$. By the fact that
$B_G(u_{w_1},1)=\{v_0,w_1\}$, we have $w_1\in N_{G}(v_0)$, which contradicts
with $|N_{G}(v_0)|=\{w_0\}$. Thus, we have $u_w=w$ for all
$w\in V_G\setminus\{v_0\}$ and then $N_{G}(w)=\{v_0\}$. This completes the
proof of part (ii). By parts (i) and (ii), we have
$$\Big(1+\frac{n-1}{n^{q(1-\alpha)}}\Big)^{{1}/{q}}\leq\|\widetilde{M}_{\alpha,G}\|_{L^p\rightarrow L^q}\leq\Big(1+\frac{n-1}{2^{q(1-\alpha)}}\Big)^{{1}/{q}}.$$
\end{proof}

Before presenting the proof of Theorem \ref{thm1.2}, let us present the
following lemma, which can be proved by the arguments similar to those
used to derive \cite[Lemma 3.1]{XZ}. The details are omitted.
	
\begin{lemma}\label{lem2.1}
Let $G=(V_G,E_G)$ be a finite connected graph with $n\ge 2$ vertices,
and let $0\leq\alpha<1$ and $0 <p<+\infty$. There exists a function
$f:V_G\rightarrow\mathbb{R}^+$ such that
$$\|\widetilde{M}_{\alpha,G}\|_{L^p\rightarrow L^p}=\frac{\|\widetilde{M}_{\alpha,G}f\|_p}{\|f\|_p}.$$
\end{lemma}

Now we present the proof of Theorem \ref{thm1.2}.
\begin{proof}[{\rm Proof of Theorem \ref{thm1.2}}]
Let $S_n=(V_{S_n},E_{S_n})$ with $V_{S_n}=\{1,\ldots,n\}$ and
$E_{S_n}=\{1\sim i:i\in\{2,\ldots,n\}\}$. The proof will be divided into
three steps:

{\bf Step 1.} At first we prove that there exists a function
$f:V_{S_n}\rightarrow\mathbb{R}^{+}$ such that (1.2) holds and
$f(1)=\max_{2\leq i\leq n}f(i)$.

In view of Lemma \ref{lem2.1}, there exists a function $f\geq0$ such that
$$\|\widetilde{M}_{\alpha,G}\|_{L^p\rightarrow L^p}=\frac{\|\widetilde{M}_{\alpha,G}f\|_p}{\|f\|_p}.\eqno(2.3)$$

If $f(1)=\max_{2\leq i\leq n}f(i)$, there is nothing to do.
If not, without loss of generality we may assume
$$f(2)\geq f(3)\geq\ldots \geq f(k)\geq f(1)\geq f(k+1)\geq\ldots\geq f(n).$$
We define the following function
\[
F(j) =
\begin{cases}
f(2), & j = 1, \\
f(1), & j = 2, \\
f(j), & j\in\{3,\ldots,n\}.
\end{cases}
\]
Clearly, $F(1) \geq F(j)$ for all $j\in\{2,3,\dots,n\}$. We will prove that
$$\|\widetilde{M}_{\alpha,S_n}\|_{L^p\rightarrow L^p}=\frac{\|\widetilde{M}_{\alpha,S_n} F\|_p}{\|F\|_p}.\eqno(2.4)$$
For any function $g:V_{S_n}\rightarrow\mathbb{R}^+$, we set
$$m_\alpha(g)=n^{\alpha-1}\sum\limits_{i=1}^ng(i),\ \ \  \ C_\alpha(g):=\max\Big\{\frac{g(2)+g(1)}{2^{1-\alpha}}, m_\alpha(g)\Big\}.$$
By the definition of $\widetilde{M}_{\alpha,S_n}$, we can write
$$\begin{array}{ll}
&\big(\widetilde{M}_{\alpha,S_n}F(1)\big)^p+\big(\widetilde{M}_{\alpha,S_n}F(2)\big)^p
=\max\big\{F(1)^p,C_\alpha(F)^p\big\}+\max\big\{F(2)^p,C_\alpha(F)^p\big\}\\
&\qquad\qquad\qquad\qquad\qquad\qquad\quad=\max\big\{f(2)^p, C_\alpha(f)^p\big\} + \max\big\{f(1)^p,C_\alpha(f)^p\big\},
\end{array}$$
$$\big(\widetilde{M}_{\alpha,S_n}f(1)\big)^p + \big(\widetilde{M}_{\alpha,S_n}f(2)\big)^p
=\max\big\{f(1)^p, C_\alpha(f)^p\big\}+ \max\big\{f(2)^p,C_\alpha(f)^p\big\},$$
$$\begin{array}{ll}
&\widetilde{M}_{\alpha,S_n}F(j)=\displaystyle\max\Big\{F(j),\frac{F(j)+F(1)}{2^{1-\alpha}},\, m_\alpha(F)\Big\}\\
&\qquad\qquad\quad=\displaystyle\max\Big\{f(j),\frac{f(j)+f(2)}{2^{1-\alpha}},\, m_\alpha(f)\Big\}\geq \widetilde{M}_{\alpha,S_n}f(j),\ \ \ j\in\{3,\ldots,n\}.
\end{array}$$
It follows that
$$\sum\limits_{j=1}^n \big(\widetilde{M}_{\alpha,S_n}F(j)\big)^p\geq\sum\limits_{j=1}^n\big(\widetilde{M}_{\alpha,S_n}f(j)\big)^p.$$
This yields that
$$\|\widetilde{M}_{\alpha,S_n}\|_{L^p\rightarrow L^p}\geq\frac{\|\widetilde{M}_{\alpha,S_n}F\|_p}{\|F\|_p}\geq\frac{\|\widetilde{M}_{\alpha,S_n}f\|_p}{\|f\|_p}
=\|\widetilde{M}_{\alpha,S_n}\|_{L^p\rightarrow L^p}.$$
So (2.4) holds and {\bf Step 1} is completed.

In what follows, without loss of generality we may assume
$$f(1) \geq f(2) \geq \ldots \geq f(n-1) \geq f(n).$$

{\bf Step 2.} We shall prove that $f|_{V\setminus\{1\}}$ takes at most
three values.

Since $f(1)\geq f(n)$ and $m_\alpha(f)\geq f(n)$, it is enough to consider
the following ten cases.
	
{\it Case 1.} ($f(n)\leq f(2)\leq f(1)\leq 2^{1-\alpha}m_\alpha(f)-f(1)\leq
m_\alpha(f)$ or $f(n)\leq f(2)\leq 2^{1-\alpha}m_\alpha(f)-f(1)\leq f(1)\leq
m_\alpha(f)$).

In this case, we get
$$\widetilde{M}_{\alpha,S_n}f(j)=m_\alpha(f), \ \ \ j\in\{1,\ldots,n\}.$$
Taking $F(i)=\frac{1}{n}\sum_{j=1}^nf(j)$ for all $i\in V_{S_n}$.
By the arguments same as case 1 in the proof of \cite[Theorem 3.2]{XZ},
we have that (1.2) holds for  $\tilde{f}$.

{\it Case 2.} ($f(n)\leq 2^{1-\alpha}m_\alpha(f)-f(1)\leq f(2)\leq f(1)
\leq m_\alpha(f)$).

In this case we may assume that there exists $k\in\{2,\dots,n-1\}$
such that
$$m_\alpha(f)\geq f(1)\geq\cdots\geq f(k)\geq 2^{1-\alpha}m_\alpha(f)-f(1)\geq f(k+1)\geq\cdots\geq f(n).$$
Thus we have
\[
\widetilde{M}_{\alpha,S_n}f(j)= \left\{
\begin{array}{ll}
2^{\alpha-1}(f(2)+f(1)),&\ j=1;\\
2^{\alpha-1}(f(j)+f(1)),&\ j\in\{2,\ldots,k\};\\
m_\alpha(f),&\ j\in\{k+1,\ldots,n\}.
\end{array}
\right.\eqno(2.5)
\]
By Jensen's inequality and (2.3), we have
$$\begin{array}{ll}
&\displaystyle(n-k)m_\alpha(f)^p+\sum\limits_{j=2}^k(2^{\alpha-1}(f(j)+f(1)))^p+(2^{\alpha-1}(f(2)+f(1)))^p\\
&=\displaystyle\|\widetilde{M}_{\alpha,S_n}\|_{L^p\rightarrow L^p}^p\sum\limits_{j=1}^n f(j)^p \\
&\geq\displaystyle\|\widetilde{M}_{\alpha,S_n}\|_{L^p\rightarrow L^p}^p\Big(\sum\limits_{j=1}^kf(j)^p+(n-k)\Big(\frac{\sum_{i=k+1}^n f(i)}{n-k}\Big)^p\Big).
\end{array}$$
Take the function
\[
\tilde{f}(j)= \left\{
\begin{array}{ll}
f(j),&\ j\in\{1,\ldots,k\};\\
\displaystyle\frac{1}{n-k}\sum\limits_{i=k+1}^n f(i),&\ j\in\{k+1,\ldots,n\}.
\end{array}
\right.
\]
Then $m_\alpha(\tilde{f})=m_\alpha(f)$, and
\[
\widetilde{M}_{\alpha,S_n}\tilde{f}(j)= \left\{
\begin{array}{ll}
2^{\alpha-1}(\tilde{f}(2)+\tilde{f}(1)),&\ j=1;\\
2^{\alpha-1}(\tilde{f}(j)+\tilde{f}(1)),&\ j\in\{2,\ldots,k\};\\
m_\alpha(\tilde{f}),&\ j\in\{k+1,\ldots,n\}.
\end{array}
\right.\eqno(2.6)
\]
It follows from (2.6) that
$$\begin{array}{ll}
&\displaystyle\sum\limits_{j=1}^n(\widetilde{M}_{\alpha,S_n}\tilde{f}(j))^p=\displaystyle
(n-k)m_\alpha(\tilde{f})^p+\sum\limits_{j=2}^k\Big(\frac{\tilde{f}(j)+\tilde{f}(1)}{2^{1-\alpha}}\Big)^p+\Big(\frac{\tilde{f}(2)+\tilde{f}(1)}{2^{1-\alpha}}\Big)^p\\
&\qquad\qquad\qquad\quad\geq\displaystyle\|\widetilde{M}_{\alpha,S_n}\|_{L^p\rightarrow L^p}^p\sum\limits_{j=1}^n \tilde{f}(j)^p,
\end{array}$$
which implies that
$$\frac{\|\widetilde{M}_{\alpha,S_n}\tilde{f}\|_p}{\|\tilde{f}\|_p}\geq\|\widetilde{M}_{\alpha,S_n}\|_{L^p\rightarrow L^p}.$$
Take
\[
g(j)= \left\{
\begin{array}{ll}
\displaystyle\Big(\frac{1}{k-2}\Big)^{1/p}\Big(\sum\limits_{i=3}^k \tilde{f}(i)^p\Big)^{{1}/{p}},&\ j\in\{3,\ldots,k\};\\
\tilde{f}(j),&\ j\in\{1,2,k+1,k+2,\ldots,n\}.
\end{array}
\right.
\]
By Jensen's inequality, one gets
$$\sum\limits_{j=1}^ng(j)=\tilde{f}(1)+\tilde{f}(2)+\sum\limits_{j=k+1}^n\tilde{f}(j)
+(k-2)\Big(\frac{1}{k-2}\Big)^{1/p}\Big(\sum\limits_{i=3}^k \tilde{f}(i)^p\Big)^{{1}/{p}}\geq\sum\limits_{j=1}^n\tilde{f}(j),$$
which leads to $m_\alpha(g)\geq m_\alpha(\tilde{f})$.

Observe that
$$\begin{array}{ll}
&\widetilde{M}_{\alpha,S_n}g(j)=\max\{g(j),2^{\alpha-1}(g(j)+g(1)),m_\alpha(g)\}\\
&\qquad\qquad=\max\{\tilde{f}(j),2^{\alpha-1}(\tilde{f}(j)+\tilde{f}(1)),m_\alpha(g)\}\geq \widetilde{M}_{\alpha,S_n}\tilde{f}(j),\ \ j\in\{k+1,\ldots,n\},
\end{array}$$
$$\begin{array}{ll}
&\widetilde{M}_{\alpha,S_n}g(j)=\max\{g(j),2^{\alpha-1}(g(j)+g(1)),m_\alpha(g)\}\\
&\qquad\qquad\geq 2^{\alpha-1}(g(j)+g(1))=\displaystyle2^{\alpha-1}\Big(\Big(\frac{1}{k-2}\Big)^{1/p}\Big(\sum\limits_{i=3}^k \tilde{f}(i)^p\Big)^{{1}/{p}}+\tilde{f}(1)\Big),\ \ j\in\{3,\ldots,k\},
\end{array}$$
$$\widetilde{M}_{\alpha,S_n}g(j)=2^{\alpha-1}(\tilde{f}(1)+\tilde{f}(2))=\widetilde{M}_{\alpha,S_n}\tilde{f}(j),\ \ j\in\{1,2\}.$$
By Minkowski inequality and (2.6), we obtain
$$\begin{array}{ll}
&\displaystyle\Big(\sum\limits_{j=3}^k(\widetilde{M}_{\alpha,S_n}g(j))^p\Big)^{1/p}\geq(k-2)^{1/p}2^{\alpha-1}\Big((k-2)^{-1/p}\Big(\sum\limits_{i=3}^k\tilde{f}(i)^p\Big)^{1/p}+\tilde{f}(1)\Big)\\
&\qquad\qquad\qquad\qquad\quad\geq\displaystyle\Big(\sum\limits_{j=3}^k\Big(\frac{\tilde{f}(j)+\tilde{f}(1)}{2^{1-\alpha}}\Big)^p\Big)^{1/p}=\Big(\sum\limits_{j=3}^k(\widetilde{M}_{\alpha,S_n}\tilde{f}(j))^p\Big)^{1/p}.
\end{array}$$
Clearly, $\|g\|_p=\|\tilde{f}\|_p$. Thus, we conclude that
$$\frac{\|\widetilde{M}_{\alpha,S_n}g\|_p}{\|g\|_p}\geq\frac{\|\widetilde{M}_{\alpha,S_n}\tilde{f}\|_p}{\|\tilde{f}\|_p}\geq\|\widetilde{M}_{\alpha,S_n}\|_{\ell^p\rightarrow\ell^p}.$$
In this case we may assume that $f(i)=f(3)$ for all $i\in\{3,\ldots,k\}$
and $f(j)=f(n)$ for all $j\in\{k+1,\ldots,n\}$.
	
{\it Case 3.} ($2^{1-\alpha}m_\alpha(f)-f(1)\leq f(n)\leq f(2)\leq f(1)\leq m_\alpha(f)$).

In this case we can write
\[
\widetilde{M}_{\alpha,S_n}f(j)= \left\{
\begin{array}{ll}
2^{\alpha-1}(f(2)+f(1)),&\ j=1;\\
2^{\alpha-1}(f(j)+f(1)),&\ j\in\{2,\ldots,n\}.
\end{array}
\right.
\]	
Taking the function
\[
\tilde{f}(j)= \left\{
\begin{array}{ll}
f(j),&\ j\in\{1,2\};\\
\displaystyle(n-2)^{-1/p}\Big(\sum\limits_{i=3}^n f(i)^p\Big)^{{1}/{p}},&\ j\in\{3,\ldots,n\}.
\end{array}
\right.
\]	
By the arguments similar to Case 2, we have
$$\frac{\|\widetilde{M}_{\alpha,S_n}\tilde{f}\|_p}{\|\tilde{f}\|_p}\geq\frac{\|\widetilde{M}_{\alpha,S_n}f\|_p}{\|f\|_p}=\|\widetilde{M}_{\alpha,S_n}\|_{L^p\rightarrow L^p}.$$
Therefore, in this case we may assume that $f(i)=f(3)$ for all
$i\in\{3,\ldots,n\}$.
	
{\it Case 4.} ($f(n)\leq m_\alpha(f)\leq 2^{1-\alpha}m_\alpha(f)-f(1)\leq f(2)\leq f(1)$).

In this case we may assume that there exists $k\in\{2,\dots,n-1\}$ such that
$$f(1)\geq \cdots \geq f(k) \geq 2^{1-\alpha}m_\alpha(f)-f(1)\geq f(k+1)\geq\cdots\geq f(n).$$
We can write
\[
\widetilde{M}_{\alpha,S_n}f(j)= \left\{
\begin{array}{ll}
\max\{f(1),2^{\alpha-1}(f(2)+f(1))\},&\ j=1;\\
2^{\alpha-1}(f(j)+f(1)),&\ j\in\{2,\ldots,k\};\\
m_\alpha(f),&\ j\in\{k+1,\ldots,n\}.
\end{array}
\right.
\]	
We consider two subcases: If $f(1)>2^{\alpha-1}(f(2)+f(1))$, then
we take the function
\[
\tilde{f}(j)= \left\{
\begin{array}{ll}
f(j),&\ j=1,2;\\
\displaystyle(k-2)^{-1/p}\Big(\sum\limits_{i=3}^k f(i)^p\Big)^{{1}/{p}},&\ j\in\{3,\ldots,k\};\\
\dfrac{1}{n-k}\sum_{i=k+1}^n f(i),&\ j\in\{k+1,\ldots,n\}.
\end{array}
\right.
\]	
By the arguments same as case 4 in the proof of \cite[Theorem 3.2]{XZ},
there exists a function $f:V_{S_n}\rightarrow\mathbb{R}$ satisfying
$\|\widetilde{M}_{\alpha,S_n}\|_{L^p\rightarrow L^p}=
\frac{\|\widetilde{M}_{\alpha,S_n}\tilde{f}\|_p}{\|\tilde{f}\|_p}$.

If $f(1)\leq 2^{\alpha-1}(f(2)+f(1))$, then we can write
\[
\widetilde{M}_{\alpha,S_n}f(j)= \left\{
\begin{array}{ll}
2^{\alpha-1}(f(2)+f(1)),&\ j=1;\\
2^{\alpha-1}(f(j)+f(1)),&\ j\in\{2,\ldots,k\};\\
m_\alpha(f),&\ j\in\{k+1,\ldots,n\}.
\end{array}
\right.
\]		
This case can be dealt with Case 2.
	
{\it Case 5.} ($f(n)\leq m_\alpha(f)\leq f(2)\leq 2^{1-\alpha}m_\alpha(f)-f(1)\leq f(1)$
or $f(n)\leq f(2)\leq m_\alpha(f)\leq 2^{1-\alpha}m_\alpha(f)-f(1)\leq f(1)$).

In this case we can write
\[
\widetilde{M}_{\alpha,S_n}f(j)= \left\{
\begin{array}{ll}
f(1),&\ j=1;\\
m_\alpha(f),&\ j\in\{2,\ldots,n\}.
\end{array}
\right.
\]
This case is same as the case 7 in the proof of \cite[Theorem 3.2]{XZ}.

{\it Case 6.} ($f(n) \leq 2^{1-\alpha}m_\alpha(f)-f(1)\leq m_\alpha(f)
\leq f(2)\leq f(1)$ or $f(n)\leq 2^{1-\alpha}m_\alpha(f)-f(1)\leq f(2)
\leq m_\alpha(f)\leq f(1)$).

In this case we may assume that there exists $l\in\{2,\ldots,n-1\}$
such that
$$f(1)\geq\cdots\geq f(l)\geq 2^{1-\alpha}m_\alpha(f)-f(1)\geq f(l+1)\geq\cdots\geq f(n).$$
This case can be dealt with Case 4.
	
{\it Case 7.} ($f(n)\leq f(2)\leq 2^{1-\alpha}m_\alpha(f)-f(1)\leq
m_\alpha(f)\leq f(1)$).

In this case we can write
\[
\widetilde{M}_{\alpha,S_n}f(j)= \left\{
\begin{array}{ll}
f(1),&\ j=1;\\
m_\alpha(f),&\ j\in\{2,\ldots,n\}.
\end{array}
\right.
\]
This case can be dealt with Case 5.

{\it Case 8.} ($2^{1-\alpha}m_\alpha(f)-f(1)\leq f(n)\leq f(2)\leq
m_\alpha(f)\leq f(1)$ or $2^{1-\alpha}m_\alpha(f)-f(1)\leq f(n)
\leq m_\alpha(f)\leq f(2)\leq f(1)$).

In this case we can write
\[
\widetilde{M}_{\alpha,S_n}f(j)= \left\{
\begin{array}{ll}
\max\{f(1), 2^{\alpha-1}(f(2)+f(1))\},&\ j=1;\\
2^{\alpha-1}(f(j)+f(1)),&\ j\in\{2,\ldots,n\}.
\end{array}
\right.
\]

We consider two subcases: If $f(1)>2^{\alpha-1}(f(2)+f(1))$, then
we take the function
\[
\tilde{f}(j)= \left\{
\begin{array}{ll}
f(j),&\ j\in\{1,2\};\\
\displaystyle(n-2)^{-1/p}\Big(\sum\limits_{i=3}^k f(i)^p\Big)^{{1}/{p}},&\ j\in\{3,\ldots,n\}.
\end{array}
\right.
\]	
By the arguments same as case 4 in the proof of \cite[Theorem 3.2]{XZ},
there exists a function $f:V_{S_n}\rightarrow\mathbb{R}$ satisfying
$\|\widetilde{M}_{\alpha,S_n}\|_{L^p\rightarrow L^p}=
\frac{\|\widetilde{M}_{\alpha,S_n}\tilde{f}\|_p}{\|\tilde{f}\|_p}$.
If $f(1)\leq 2^{\alpha-1}(f(2)+f(1))$, then we can write
\[
\widetilde{M}_{\alpha,S_n}f(j)= \left\{
\begin{array}{ll}
2^{\alpha-1}(f(2)+f(1)),&\ j=1;\\
2^{\alpha-1}(f(j)+f(1)),&\ j\in\{2,\ldots,n\}.
\end{array}
\right.
\]		
This case is similar to Case 2.

{\it Case 9.} ($f(n)\leq f(2)\leq m_\alpha(f)\leq f(1)\leq 2^{1-\alpha}m_\alpha(f)-f(1)$
or $f(n)\leq m_\alpha(f)\leq f(2)\leq f(1)\leq 2^{1-\alpha}m_\alpha(f)-f(1)$).

In this case we can write
\[
\widetilde{M}_{\alpha,S_n}f(j)= \left\{
\begin{array}{ll}
f(1),&\ j=1;\\
m_\alpha(f),&\ j\in\{2,\ldots,n\}.
\end{array}
\right.
\]		
This case can be dealt with Case 5.

{\it Case 10.} ($f(n)\leq f(1)\leq m_\alpha(f)\leq 2^{1-\alpha}m_\alpha(f)-f(1)$).

In this case we have $\widetilde{M}_{\alpha,S_n}f(j)=m_\alpha(f)$ for
$1\leq j\leq n$. Set $\tilde{f}(i)=\frac{1}{n}\sum_{j=1}^nf(j)$ for
$1\leq i\leq n$. As in the proof of Case 1, we can conclude that
$$\frac{\|\widetilde{M}_{\alpha,S_n}\tilde{f}\|_p}{\|\tilde{f}\|_p}=\|\widetilde{M}_{\alpha,S_n}\|_{L^p\rightarrow L^p}.$$
This completes the proof of Theorem \ref{thm1.2}.
\end{proof}
	
\bigskip

\section{Sharp $p$-variation constants for $M_{S_n}$ and $\widetilde{M}_{S_n}$
--Proofs of Theorems \ref{thm1.3} and \ref{thm1.4}}\label{S3}

\begin{lemma}\label{lem3.1}
Let $0<p<1$ and $r>\frac{1}{2}$. Define the following function
$$\phi(t)=\min\{2^{-1}t,1\}-r^pt^p,\ \ \ t\in[0,\infty).$$
Then for any $u,\,v>0$, there exist $a,\,b$ such that
$$\phi(u)+\phi(v)\leq\phi(a)+\phi(b),\ \ \ \{a,b\}\cap\{0,2\}\neq\emptyset,\ \ \ a+b=u+v.\eqno(3.1)$$
\end{lemma}
\begin{proof}
Let us fix $u,\,v>0$ and we define the following function
$$G(t)=\phi(t)+\phi(u+v-t),\ \ \ t\in[0,u+v].$$

We consider two cases:

{\it Case 1.} ($u+v<4$). In this case we may assume
that $u<2$. By the definition of $\phi$, we see that
$G(t)$ is convex on $[\max\{0,u+v-2\},\min\{2,u+v\}]$.
If $u+v\leq2$, then we have $[\max\{0,u+v-2\},
\min\{2,u+v\}]=[0,u+v]$ and
$$G(u)\leq\max\{G(0),G(u+v)\}=G(0).$$
In this case we take $a=0$ and $b=u+v$.
If $u+v\in(2,4)$, then we have $[\max\{0,u+v-2\},
\min\{2,u+v\}]=[u+v-2,2]$ and
$$G(u)\leq\max\{G(u+v-2),G(2)\}=G(2).$$
In this case we taking $a=2$ and $b=u+v-2$. This
prove (5.1) in this case.

{\it Case 2.} ($u+v\geq4$). If $u=2$ or $v=2$, then there is nothing
to do. If not, without loss of generality we may assume $v>2$. If
$u>2$, then $G(t)$ is convex on $[2,u+v-2]$ and $v\in[2,u+v-2]$. So
$$G(v)\leq\max\{G(2),G(u+v-2)\}=G(2)=\phi(2)+\phi(u+v-2).$$
In this case we take $a=2$ and $b=u+v-2$. If $u<2$,
then we set $\theta=\frac{u}{2}\in(0,1)$. Clearly,
$$v=(1-\theta)(u+v)+\theta(u+v-2),\ \ u=(1-\theta)0+2\theta.$$
Since $\phi(t)$ is convex on $[0,2]$. Moreover,
$\phi(t)$ is also convex on $[0,2]$. Thus,
$$\begin{array}{ll}
&\phi(u)+\phi(v)\leq(1-\theta)\phi(0)+\theta\phi(2)
+(1-\theta)\phi(u+v)+\theta\phi(u+v-2)\\
&\qquad\qquad\quad\leq(1-\theta)\phi(u+v)+\theta(\phi(2)+\phi(u+v-2))\\
&\qquad\qquad\quad\leq\max\{\phi(u+v),\phi(2)+\phi(u+v-2)\}.
\end{array}$$
If $\phi(u+v)\leq \phi(2)+\phi(u+v-2)$, then we take $a=2$ and $b=u+v-2$.
If $\phi(u+v)>\phi(2)+\phi(u+v-2)$, then we take $a=0$ and $b=u+v$. This
proves (3.1) in this case.
\end{proof}

We now present the proofs of Theorems \ref{thm1.3} and \ref{thm1.4}.

\begin{proof}[{\rm Proof of Theorem \ref{thm1.3}}]
Without loss of generality we may assume $S_n=(V_{S_n},E_{S_n})$, where
$V_{S_n}=\{1,2,\ldots,n\}$ and $E_{S_n}=\{1\sim 2,\ldots,1\sim n\}$.
Clearly,
$${\rm Var}_{p,S_n}(\delta_2)=1,\ \ \ M_{S_n}\delta_2(2)=1,\ \ \ M_{S_n}\delta_2(i)=\frac{1}{n},\ \ \ \ i\in\{1,3,4,\ldots,n\}.$$
It is clear that
$${\bf C}_{S_n}^{p,p}\geq\frac{{\rm Var}_{p,S_n}(M_{S_n}\delta_2)}{{\rm Var}_{p,S_n}(\delta_2)}=1-\frac{1}{n}.$$
Therefore, for (1.3) it remains to show that for any $f:V_{S_n}
\rightarrow\mathbb{R}$ satisfying ${\rm Var}_{p,S_n}(f)>0$,
$${\rm Var}_{p,S_n}(M_{S_n}f)\leq\Big(1-\frac{1}{n}\Big){\rm Var}_{p,S_n}(f).\eqno(3.1)$$

Let $f:V_{S_n}\rightarrow\mathbb{R}$ satisfying ${\rm Var}_{p,S_n}(f)>0$.
Without loss of generality we may assume $f\geq0$. For convenience, we set
$$a_i=f(i),\ \  d_i=a_i-a_1, \ \ \ m=\frac{1}{n}\sum\limits_{i=1}^na_i, \ \ \ i\in\{1,\ldots,n\}.$$
It is clear that
$$M_{S_n}f(1)=\max\{a_1,m\},\ \ \ Mf(i)=\max\Big\{a_i,\frac{a_i+a_1}{2},m\Big\},\ \ \ \ i\in\{2,\ldots,n\}.$$

Next we consider two steps:

{\bf Step 1.} Proof of (3.1) for the case $a_1\leq m$.

Without loss of generality we may assume $a_i>m$ for $i\in\{k+1,
\ldots,n\}$ and $a_j\leq m$ for $j\in\{1,\ldots,k\}$. Clearly,
$k\in[1,n-2]$. Set
$$u=m-a_1,\ \ \ \ A=\sum\limits_{i=k+1}^nx_i,\ \ \ \ x_j=a_j-m,\ \ \ j\in\{k+1,\ldots,n\}.$$
Clearly, $({\rm Var}_{p,S_n}(M_{S_n}f))^p=\sum_{j=k+1}^nx_j^p$,
$nu=\sum_{j=2}^nd_j$ and $d_j=x_j+u>0$ for $j\in\{k+1,\ldots,n\}$.
Moreover,
$$\sum\limits_{j=1}^kd_j=ku-A.$$
Observe that
$$({\rm Var}_{p,S_n}(f))^p=\sum\limits_{j=1}^k|d_j|^p+\sum\limits_{j=k+1}^n(x_j+u)^p\geq F(u):=|ku-A|^p+\sum\limits_{j=k+1}^n(x_j+u)^p,$$
$$({\rm Var}_{p,S_n}(f))^p\geq F(k^{-1}A),\ \ \ u\in[k^{-1}A,\infty).$$
We also note that $F(u)$ is concave on $u\in[1,k^{-1}A]$. Therefore, we
have
$$({\rm Var}_{p,S_n}(f))^p\geq \min\{F(k^{-1}A),F(0)\}.$$
Note that
$$F(k^{-1}A)=\sum\limits_{j=k+1}^n\Big(x_j+\frac{A}{k}\Big)^p
\geq\displaystyle\sum\limits_{j=k+1}^nx_j^p\Big(1+\frac{1}{k}\Big)^p\geq\Big(\frac{n}{n-1}\Big)^p({\rm Var}_{p,S_n}(M_{S_n}f))^p,$$
$$F(0)=|A|^p+\sum\limits_{j=k+1}^n(x_j)^p\geq\Big(1+\frac{1}{k}\Big)\sum\limits_{j=k+1}^nx_j^p\geq\Big(\frac{n}{n-1}\Big)^p({\rm Var}_{p,S_n}(M_{S_n}f))^p.$$
So
$$({\rm Var}_{p,S_n}(f))^p\geq \Big(\frac{n}{n-1}\Big)^p({\rm Var}_{p,S_n}(M_{S_n}f))^p.$$
This yields (3.1) in this case.

{\bf Step 2.} Proof of (3.1) for the case $a_1>m$.

We set $b=a_1-m$. Denote
$$I^{+}:=\{i\in\{2,\ldots,n\}:a_i>a_1\},\ \ \ \ I^{-}:=\{i\in\{2,\ldots,n\}:a_i<a_1\}.$$
Set $x_i=d_i=a_i-a_1$ for $i\in I^{+}$ and $y_j=-d_j=a_1-a_j$ for $j\in I^{-}$.
Clearly, $\sum_{j\in I^{-}}y_j-\sum_{i\in I^{+}}x_i=nb$ and
$$({\rm Var}_{p,S_n}(M_{S_n}f))^p=\sum\limits_{i\in I^{+}}x_i^p+\sum\limits_{j\in I^{-}}\min\Big\{b,\frac{y_j}{2}\Big\}^p.$$
For convenience, we set $k=|I^{+}|$ and $\ell=|I^{-}|$. If $k=0$, then
we have
$$({\rm Var}_{p,S_n}(M_{S_n}f))^p=\sum\limits_{j\in I^{-}}\min\Big\{b,\frac{y_j}{2}\Big\}^p\leq\frac{1}{2^p}\sum\limits_{j\in I^{-}}y_j^p
=\frac{1}{2^p}({\rm Var}_{p,S_n}(f))^p\leq\Big(1-\frac{1}{n}\Big)^p({\rm Var}_{p,S_n}(f))^p.$$
There is nothing to do. In what follows, we assume $1\leq k\leq n-\ell-1$.
Set $\tilde{x_i}=\frac{x_i}{b}$ and $\tilde{y_j}=\frac{y_j}{b}$.
Clearly, $\tilde{x_i}>0$ for any $i\in I^{+}$ and $y_j>0$ for any
$j\in I^{-}$. Moreover, $\sum_{j\in I^{-}}\tilde{y_j}-\sum_{i\in I^{+}}\tilde{x_i}=n$,
$({\rm Var}_{p,S_n}(f))^p=b^p(\sum_{i\in I^{+}}\tilde{x_i}^p+
\sum_{j\in I^{-}}\tilde{y_j}^p)$ and
$$({\rm Var}_{p,S_n}(M_{S_n}f))^p=b^p\Big(\sum\limits_{i\in I^{+}}\tilde{x_i}^p+\sum\limits_{j\in I^{-}}\min\Big\{1,\frac{\tilde{y_j}}{2}\Big\}^p\Big).$$
Let
$$\phi(t)=\min\Big\{\frac{t}{2},1\Big\}-\Big(\frac{n-1}{n}\Big)^pt^p,\ \ \ t\in[0,\infty).$$
It suffices to show that
$$\Big(1-\Big(\frac{n-1}{n}\Big)^p\Big)\sum\limits_{i\in I^{+}}\tilde{x_i}^p+\sum\limits_{j\in I^{-}}\phi_n(\tilde{y_j})\leq0,\ \ \ \sum_{j\in I^{-}}\tilde{y_j}=\sum_{i\in I^{+}}\tilde{x_i}+n.\eqno(3.2)$$

Set $T=\sum_{i\in I^{+}}\tilde{x_i}$, then $\sum_{j\in I^{-}}\tilde{y_j}=n+T>n$.
Next we shall prove that there exists $r\in\{0,1,\ldots,\ell-1\}$ such that
$n+T-2r\geq0$ and
$$\sum\limits_{j\in I^{-}}\phi(\tilde{y_j})\leq r\phi(2)+\phi(T+n-2r).\eqno(3.3)$$
Without loss of generality we may assume that
$I^{-}=\{i_1,\ldots,i_{\ell}\}$ with $i_1<i_2<\ldots<i_{\ell}$
and
$$\tilde{y_{i_1}}\leq \tilde{y_{i_2}}\leq \tilde{y_{i_3}}\leq\cdots\leq\cdots\leq \tilde{y_{i_{\ell}}}.$$
Note that $\frac{n-1}{n}>\frac{1}{2}$ for $n\geq3$. By Lemma
\ref{lem3.1}, there exist $\{z_j\}_{j=1}^{\ell}$ such that
$\{z_1,\ldots,z_{\ell-1}\}\subset\{0,2\}$ and $z_\ell=T+n-2r$
with $r=|\{i:1\leq i\leq\ell-1:z_j=2\}|$ and
$$\sum\limits_{j\in I^{-}}\phi(\tilde{y_j})\leq
\sum\limits_{j\in I^{-}}\phi(z_j)\leq r\phi(2)+\phi(T+n-2r).$$
This yields (3.2).

Since $r\leq\ell$, then we have $k\leq n-1-r$. By Jensen's
inequality, one gets
$$\sum\limits_{j\in I^{+}}x_j^p\leq k^{1-p}T^p\leq(n-1-r)^{1-p}T^p.$$
In view of (3.3), for (3.2) it suffices to show that
$$\Phi(T,r):=\Big(1-\Big(\frac{n-1}{n}\Big)^p\Big)(n-2-r)^{1-p}T^p-
r\Big(\Big(\frac{2(n-1)}{n}\Big)^p-1\Big)+\phi(T+n-2r)\leq0,\eqno(3.4)$$
where $r\in\{0,\ldots,n-3\}$ and $T>\{0,2r-n\}$.

We consider three cases:

{\it Case 1.} ($\frac{n}{2}\leq r\leq n-3$). At first we conclude that
$$\Phi(T,r)\leq\max\{\Phi(2r-n,r),\Phi(2+2r-n,r)\}.\eqno(3.5)$$
If $T+n-2r\geq2$, then we can write
$$\Phi(T,r):=\Big(1-\Big(\frac{n-1}{n}\Big)^p\Big)(n-2-r)^{1-p}T^p-
r\Big(\Big(\frac{2(n-1)}{n}\Big)^p-1\Big)+1-\Big(\frac{n-1}{n}\Big)^p(T+n-2r)^p.$$
It is clear that
$$\begin{array}{ll}
&\frac{\partial\Phi(T,r)}{\partial T}=\displaystyle pT^{p-1}\Big(1-\Big(\frac{n-1}{n}\Big)^p\Big)(n-2-r)^{1-p}
-p(T+n-2r)^{p-1}\Big(\frac{n-1}{n}\Big)^p\\
&\qquad\quad\leq\displaystyle pT^{p-1}
\Big(\Big(1-\Big(\frac{n-1}{n}\Big)^p\Big)(n-2-r)^{1-p}-\Big(\frac{n-1}{n}\Big)^p\Big)\\
&\qquad\quad\leq\displaystyle pT^{p-1}\Big(\frac{n-1}{n}\Big)^p
\Big(\Big(\Big(\frac{n}{n-1}\Big)^p-1\Big)(n-2)^{1-p}-1\Big).
\end{array}$$
By the differentiable mean value theorem, it is not difficult to see that
$$0<\Big(\frac{n}{n-1}\Big)^p-1\leq\frac{p}{n-1}.$$
It follows that
$$\Big(\Big(\frac{n}{n-1}\Big)^p-1\Big)(n-2)^{1-p}\leq\frac{p(n-2)}{n-1}<1.$$
So $\frac{\partial\Phi(T,r)}{\partial T}<0$. This yields (3.5) in this case.

If $T+n-2r\in[0,2]$, then we can write
$$\begin{array}{ll}
&\Phi(T,r):=\displaystyle\Big(1-\Big(\frac{n-1}{n}\Big)^p\Big)(n-2-r)^{1-p}T^p\\
&\qquad\qquad\quad-\displaystyle r\Big(\Big(\frac{2(n-1)}{n}\Big)^p-1\Big)+\Big(\Big(\frac{n-1}{n}\Big)^p-2^{-p}\Big)(T+n-2r)^p.
\end{array}$$
Clearly,
$$\frac{\partial\Phi(T,r)}{\partial T}=pT^{p-1}g(T),$$
where
$$g(T):=\Big(\Big(1-\Big(\frac{n-1}{n}\Big)^p\Big)(n-2-r)^{1-p}
+\Big(\Big(\frac{n-1}{n}\Big)^p-2^{-p}\Big)\Big(1+\frac{n-2r}{T}\Big)^{p-1}\Big).$$
Since $g(T)$ is increasing on $[2r-n,2+2r-n]$ and $pT^{p-1}>0$, then (3.5)
holds in this case.

Next we prove that
$$\max\{\Phi(2r-n,r),\Phi(2+2r-n,r)\}\leq0.\eqno(3.6)$$
By the differentiable mean value theorem, one gets
$$\Big(1-\Big(\frac{n-1}{n}\Big)^p\Big)-\Big(\Big(\frac{2(n-1)}{n}\Big)^p-1\Big)
\leq\displaystyle\frac{p}{n-1}-p\ln\frac{2(n-1)}{n}\leq \frac{p}{n-1}-\frac{p(n-2)}{2(n-1)}\leq0.
\eqno(3.7)$$
In view of (3.7), we have
$$\Phi(2r-n,r)\leq r\Big(\Big(1-\Big(\frac{n-1}{n}\Big)^p\Big)-
\Big(\Big(\frac{2(n-1)}{n}\Big)^p-1\Big)\Big)\leq 0,$$
$$\begin{array}{ll}
&\Phi(2+2r-n,r)\leq\displaystyle\Big(1-\Big(\frac{n-1}{n}\Big)^p\Big)(n-2-r)^{1-p}(2+2r-n)^p-
r\Big(\Big(\frac{2(n-1)}{n}\Big)^p-1\Big)+\phi(2)\\
&\qquad\qquad\qquad\quad\leq\displaystyle(r+1)\Big(\Big(1-\Big(\frac{n-1}{n}\Big)^p\Big)-
\Big(\Big(\frac{2(n-1)}{n}\Big)^p-1\Big)\Big)\leq0.
\end{array}$$
Thus, (3.7) holds.

{\it Case 2.} ($r=\frac{n-1}{2}$). In this case $n\geq5$ and $n$
is odd. Write
$$\Phi(T,r):=\Big(1-\Big(\frac{n-1}{n}\Big)^p\Big)\Big(\frac{n-3}{2}\Big)^{1-p}T^p-
\frac{n-1}{2}\Big(\Big(\frac{2(n-1)}{n}\Big)^p-1\Big)+\phi(T+1).$$
If $0<T\leq1$, then $T\leq\frac{n-3}{2}$. By (3.7) and the fact that
$\phi(t)\leq0$, one gets
$$\Phi(T,r)\leq\frac{n-3}{2}\Big(1-\Big(\frac{n-1}{n}\Big)^p\Big)-\Big(\Big(\frac{2(n-1)}{n}\Big)^p-1\Big)\leq0.$$
If $T>1$, then we can write
$$\Phi(T,r):=\Big(1-\Big(\frac{n-1}{n}\Big)^p\Big)\Big(\frac{n-3}{2}\Big)^{1-p}T^p-
\frac{n-1}{2}\Big(\Big(\frac{2(n-1)}{n}\Big)^p-1\Big)+1-\Big(\frac{n-1}{n}\Big)^p(T+1)^p.$$
Let $T_0$ be such that $\frac{\partial\Phi(T,r)}{\partial T}(T_0)=0$.
A simple computation gives that
$$\frac{T_0}{T_0+1}=\frac{n-3}{2}\Big(\Big(\frac{n}{n-1}\Big)^p-1\Big)^{\frac{1}{1-p}}
\leq\frac{p(n-3)}{2(n-1)}<\frac{1}{2},$$
which leads to $T_0<1$. It is easy to see that $\frac{\partial\Phi(T,r)}{\partial T}<0$
when $T>T_0$. So $\Phi(T,r)$ is decreasing on $[1,\infty)$ and then
$$\Phi(T,r)\leq \Phi(1,r)\leq0.$$

{\it Case 3.} ($0\leq r\leq\frac{n-2}{2}$). In this case
we can write
$$\Phi(T,r):=\Big(1-\Big(\frac{n-1}{n}\Big)^p\Big)(n-2-r)^{1-p}T^p-
r\Big(\Big(\frac{2(n-1)}{n}\Big)^p-1\Big)+1-\Big(\frac{n-1}{n}\Big)^p(T+n-2r)^p.$$
For convenience, we set
$$a=\Big(\frac{n-1}{n}\Big)^p,\ \ \ \ b=2^pa-1,\ \ \ \delta=\frac{1}{a}-1,\ \ \sigma=\frac{1}{1-p}.$$
Clearly, $\sigma>1$ and $\delta\in(0,1)$ since $p\in(0,1)$
and $n\geq4$. Then it suffices to show that
$$\Phi(T,r):=(1-a)(n-2-r)^{1-p}T^p+1-rb-a(T+n-2r)^p\leq0.\eqno(3.8)$$
For convenience, we set
$\tau=(n-2-r)\delta^{\sigma}$ and $T_1=\frac{\tau(n-2r)}{1-\tau}$.
Note that
$$\frac{\partial\Phi(T,r)}{\partial T}=pa\Big((n-2-r)^{1-p}\delta T^{p-1}-(T+L)^{p-1}\Big).$$
It is not difficult to check that $\frac{\partial\Phi(T,r)}{\partial T}>0$
when $T<T_1$ and $\frac{\partial\Phi(T,r)}{\partial T}<0$ when $T>T_1$.
Thus,
$$\Phi(T,r)\leq \Phi(T_1,r)=1-rb-a(n-2r)^p(1-(n-2-r)\delta^{\sigma})^{1-p},\ \ \forall T>0.$$
Thus, for (3.8) it suffices to show that
$$F(r):=1-rb-a(n-2r)^p(1-(n-2-r)\delta^{\sigma})^{1-p}\leq0.\eqno(3.9)$$
It is clear that
$$F'(r)=-b+2ap(n-2r)^{p-1}(1-(n-2-r)\delta^{\sigma})^{1-p}-a(1-p)\delta^{\sigma}(n-2r)^p
(1-(n-2-r)\delta^{\sigma})^{-p},$$
$$F''(r)=p(1-p)(n-2r)^p(1-(n-2-r)\delta^{\sigma})^{1-p}\Big(\frac{2}{n-2r}
+\frac{\delta^{\sigma}}{1-(n-2-r)\delta^{\sigma}}\Big)^2\geq0.$$
Then $F(r)$ is convex on $[0,\frac{n-2}{2}]$. Thus, for (3.9), it is
enough to conclude that
$$\max\Big\{F(0),F\Big(\frac{n-2}{2}\Big)\Big\}\leq0.\eqno(3.10)$$
To prove (3.10), let us firstly prove
$$(n-2)\delta+(n-1)^{-p}\leq1.\eqno(3.11)$$
Since $f(x)=t^p$ is concave on $(0,\infty)$, then we have
$$a=\Big(\frac{n-2}{n-1}+\frac{1}{(n-1)n}\Big)^p\geq\frac{n-2}{n-1}+\frac{1}{(n-1)n^p}.$$
It follows that
$$n-1\geq(n-2)a^{-1}+n^{-p}a^{-1}=(n-2)(1+\delta)+(n-1)^{-p},$$
which leads to (3.11). In view of (3.11) and the fact that
$\delta^\sigma\leq\delta$, one gets
$$F(0):=1-(n-1)^p(1-(n-2)\delta^{\sigma})^{1-p}\leq
1-(n-1)^p(1-(n-2)\delta)^{1-p}\leq 0.$$
Next we prove
$$F\Big(\frac{n-2}{2}\Big)\leq0.\eqno(3.12)$$
For convenience, we set $D=\frac{n-2}{2}$. Then we have
$$F(\frac{n-2}{2}):=1-bD-a2^p(1-\delta^{\sigma}D)^{1-p}.$$
If $bD\geq1$, then (3.12) is clear. If $bD<1$, then we have
$$a\geq\Big(\frac{D}{2D+1}+\frac{D+1}{2(2D+1)}\Big)\geq\frac{D+(D+1)2^{-p}}{2D+1},$$
since $p\in(0,1)$. It follows that
$$\Big(\frac{1-bD}{a2^p}\Big)^{\sigma}+D\delta^\sigma
\leq \frac{1-bD}{a2^p}+D\delta=\frac{D+(D+1)2^{-p}}{a}-2D\leq 1,$$
since $2^pa\geq1$. So $1-D\delta^\sigma\geq
\Big(\frac{1-bD}{a2^p}\Big)^{\sigma}$ and then
$$F\Big(\frac{n-2}{2}\Big)\leq1-bD-a2^p((1-\delta^{\sigma}D)^{1-p})
\leq1-bD-2^pa\Big(\frac{1-bD}{a2^p}\Big)\leq0.$$
This proves (3.12) in this case. So (3.10) holds. This completes the
proof of Theorem \ref{thm1.3}.
\end{proof}

\begin{proof}[{\rm Proof of Theorem \ref{thm1.4}}]
Without loss of generality we may assume $S_n=(V_{S_n},E_{S_n})$, where
$V_{S_n}=\{1,2,\ldots,n\}$ and $E_{S_n}=\{1\sim 2,\ldots,1\sim n\}$. It
is clear that $\widetilde{M}_{S_n}\delta_2(2)=1$, $\widetilde{M}_{S_n}
\delta_2(1)=\frac{1}{2}$ and $\widetilde{M}_{S_n}\delta_2(i)=\frac{1}{n}$
for any $i\in\{3,\ldots,n\}$. It follows that
$$\widetilde{{\bf C}}_{S_n}^{p,p}\geq\frac{{\rm Var}_{p,S_n}(\widetilde{M}_{S_n}\delta_2)}{{\rm Var}_{p,S_n}(\delta_2)}\geq\Big(2^{-p}+(n-2)\Big(\frac{n-2}{2n}\Big)^p\Big)^{1/p}.$$
Thus, for (1.4) it suffices to show that for any $f:\mathbb{R}\rightarrow\mathbb{R}$
satisfying ${\rm Var}_{p,S_n}(f)<\infty$,
$${\rm Var}_{p,S_n}(\widetilde{M}_{S_n}f)\leq \Big(2^{-p}+(n-2)\Big(\frac{n-2}{2n}\Big)^p\Big)^{1/p}{\rm Var}_{p,S_n}(f).\eqno(3.13)$$

Let us fix $f:\mathbb{R}\rightarrow\mathbb{R}$ satisfying ${\rm Var}_{p,S_n}(f)<\infty$.
Without loss of generality we may assume $f\geq0$ and $f(2)=\max_{2\leq i\leq m}f(i)$.
Set $m=\frac{1}{n}\sum_{i=1}^nf(i)$. It is clear that
\[
\widetilde{M}_{S_n}f(i)= \left\{
\begin{array}{ll}
\max\Big\{f(1),\frac{f(1)+f(2)}{2},m\Big\},&\ \ i=1;\\
\max\Big\{f(i),\frac{f(1)+f(i)}{2},m\Big\},&\ \ \ i\in\{2,\ldots,n\}.
\end{array}
\right.
\]
Let $d_i=f(i)-f(1)$ for $2\leq i\leq n$. For (3.13) it suffices to show that
$$|\widetilde{M}_{S_n}f(i)-\widetilde{M}_{S_n}f(1)|\leq\frac{1}{2}|d_i|+\frac{n-2}{2n}\sum\limits_{2\leq j\leq n,\atop j\neq i}|d_j|.\eqno(3.14)$$
In fact, in view of (3.14), we have
$$\begin{array}{ll}
&({\rm Var}_{p,S_n}(\widetilde{M}_{S_n}f))^p\leq\displaystyle\sum\limits_{i=2}^n\Big(\frac{1}{2^p}|d_i|^p+\Big(\frac{n-2}{2n}\Big)^p\sum\limits_{2\leq j\leq n,\atop j\neq i}|d_j|^p\Big)\\
&\qquad\qquad\qquad\quad=\displaystyle\Big(2^{-p}+(n-2)\Big(\frac{n-2}{2n}\Big)^p\Big)\sum\limits_{i=2}^n|d_i|^p\\
&\qquad\qquad\qquad\quad\leq\displaystyle\Big(2^{-p}+(n-2)\Big(\frac{n-2}{2n}\Big)^p\Big)({\rm Var}_{p,S_n}(f))^p,
\end{array}$$
which leads to (3.13).

Now we prove (3.14). We discuss two cases: If $f(1)\geq f(2)$, then we have
$$\widetilde{M}_{S_n}f(1)=f(1),\ \ \ \widetilde{M}_{S_n}f(2)=\frac{f(1)+f(2)}{2},\ \ \ \widetilde{M}_{S_n}f(i)=\max\Big\{\frac{f(1)+f(i)}{2},m\Big\},\ i\in\{3,\ldots,n\}.$$
It follows that
$$|\widetilde{M}_{S_n}f(1)-\widetilde{M}_{S_n}f(i)|=\widetilde{M}_{S_n}f(1)-\widetilde{M}_{S_n}f(i)\leq\widetilde{M}_{S_n}f(1)-\frac{f(1)+f(i)}{2}\leq\frac{1}{2}|d_i|,\ \ i\in\{2,\ldots,n\}.$$
This gives (3.14). If $f(1)\leq f(2)$, then we can write
$$\widetilde{M}_{S_n}f(1)=\frac{f(1)+f(2)}{2},\ \ \widetilde{M}_{S_n}f(2)=f(2),$$
$$\widetilde{M}_{S_n}f(i)=\max\Big\{f(i),\frac{f(1)+f(i)}{2},m\Big\},\ \ \ \ i\in\{3,\ldots,n\}.$$
Clearly,
$$|\widetilde{M}_{S_n}f(1)-\widetilde{M}_{S_n}f(2)|=f(2)-\frac{f(1)+f(2)}{2}=\frac{f(2)-f(1)}{2}\leq\frac{1}{2}|d_2|.$$
If $\widetilde{M}_{S_n}f(i)=f(i)$, then we have $f(i)\geq f(1)$ and
$$|\widetilde{M}_{S_n}f(1)-\widetilde{M}_{S_n}f(i)|=\Big|\frac{f(1)+f(2)}{2}-f(i)\Big|\leq\frac{1}{2}(|d_i|+f(2)-f(i))\leq\frac{1}{2}(|d_i|+|d_2|).$$
This gives (3.14) in this case. If $\widetilde{M}_{S_n}f(i)=\frac{f(1)+f(i)}{2}$,
then we have $f(i)\leq f(1)$ and
$$|\widetilde{M}_{S_n}f(1)-\widetilde{M}_{S_n}f(i)|=\Big|\frac{f(1)+f(2)}{2}-\frac{f(1)+f(i)}{2}\Big|\leq \frac{1}{2}(|d_i|+|d_2|).$$
This gives (3.14) in this case. If $\widetilde{M}_{S_n}f(i)=m$, then
we have
$$\begin{array}{ll}
&|\widetilde{M}_{S_n}f(1)-\widetilde{M}_{S_n}f(i)|=\displaystyle\Big|\frac{f(1)+f(2)}{2}-m\Big|\\
&\qquad\qquad\qquad\qquad\quad=\displaystyle\Big|\frac{n-2}{2n}d_2-\frac{1}{n}\sum\limits_{3\leq j\leq n}d_j\Big|\\
&\qquad\qquad\qquad\qquad\quad\leq\displaystyle\frac{n-2}{2n}|d_2|+\frac{1}{n}|d_i|+\frac{1}{n}\sum\limits_{3\leq j\leq n,\atop j\neq i}|d_j|\\
&\qquad\qquad\qquad\qquad\quad\leq\displaystyle\frac{1}{2}|d_i|+\frac{n-2}{2n}\sum\limits_{3\leq j\leq n,\atop j\neq i}|d_j|
\leq\displaystyle\frac{1}{2}|d_i|+\frac{n-2}{2n}\sum\limits_{2\leq j\leq n,\atop j\neq i}|d_j|,
\end{array}$$
which proves (3.14) in this case.
\end{proof}

\bigskip

\section{Sharp $q$-variation constants for fractional maximal operators
--Proofs of Theorems \ref{thm1.5}--\ref{thm1.7}}\label{S4}

In this section we present the proofs of Theorems \ref{thm1.5}--\ref{thm1.7}.

\begin{proof}[{\rm Proof of Theorem \ref{thm1.5}}] Let $0<\alpha<1$
and $q=\frac{1}{1-\alpha}$. Without loss of generality we may assume
$K_n=(V_{K_n},E_{K_n})$ with $V_{K_n}=\{1,\ldots,n\}$. It is clear that
$${\rm Var}_{K_n}(\delta_1)=n-1,\ \ \ \ {\rm Var}_{q,K_n}(M_{\alpha,K_n}\delta_1)=(n-1)^{1-\alpha}(1-n^{\alpha-1}).$$
It follows that
$${\bf C}_{\alpha,K_n}^{1,q}\geq\frac{{\rm Var}_{q,K_n}(M_{\alpha,K_n}\delta_1)}{{\rm Var}_{K_n}(\delta_1)}\geq (n-1)^{-\alpha}(1-n^{\alpha-1}).$$
Thus, for (1.5) it suffices to show that for any $f:V_{K_n}\rightarrow\mathbb{R}$
with ${\rm Var}_{K_n}(f)>0$ such that
$${\rm Var}_{q,K_n}(M_{\alpha,K_n}f)\leq (n-1)^{-\alpha}(1-n^{\alpha-1}){\rm Var}_{K_n}(f).\eqno(4.1)$$

Let $f:V_{S_n}\rightarrow\mathbb{R}$ with ${\rm Var}_{K_n}(f)>0$. Without
loss of generality we may assume $f\geq0$ since
$M_{\alpha,K_n}f=M_{\alpha,K_n}|f|$ and ${\rm Var}_{K_n}(|f|)
\leq{\rm Var}_{K_n}(f)$. We may assume that
$$f(n)\geq f(n-1)\cdots\geq f(2)\geq f(1).$$
Clearly,
$$M_{\alpha,K_n}f(i)=\max\Big\{f(i),n^{\alpha-1}\sum\limits_{j=1}^nf(j)\Big\},\ \ i\in\{1,\ldots,n\}.$$
If $f(n)\leq n^{\alpha-1}\sum_{j=1}^nf(j)$, then we have $M_{\alpha,K_n}f(i)=
n^{\alpha-1}\sum_{j=1}^nf(j)$ for all $i\in\{1,\ldots,n\}$ and (4.1)
is clear.

In what follows, we assume that
$$f(1)\leq f(2)\leq\cdots\leq f(m)\leq n^{\alpha-1}\sum\limits_{j=1}^nf(j)\leq f(m+1)\leq\cdots\leq f(n),$$
where $m\in[1,n-1]$. At first we prove
$${\rm Var}_{K_n}(M_{\alpha,K_n}f)\le(1-n^{\alpha-1}){\rm Var}_{K_n}(f).\eqno(4.2)$$
We can write
$$\begin{array}{ll}
&{\rm Var}_{K_n}(f)=\displaystyle\sum\limits_{1\leq i<j\leq n}(f(j)-f(i))\\
&\qquad\qquad=\displaystyle\sum\limits_{1\leq i<j\leq m}(f(j)-f(i))+\sum\limits_{m+1\leq i<j\leq n}(f(j)-f(i))+\sum\limits_{i=1}^{m}\sum\limits_{j=m+1}^{n}(f(j)-f(i))\\
&\qquad\qquad=\displaystyle\sum\limits_{1\leq i<j\leq m}(f(j)-f(i))+\sum\limits_{m+1\leq i<j\leq n}(f(j)-f(i))+m\sum_{j=m+1}^nf(j)-(n-m)\sum\limits_{i=1}^mf(i).
\end{array}$$
$${\rm Var}_{K_n}(M_{\alpha,K_n}f)=\sum\limits_{m+1\leq i<j\leq n}(f(j)-f(i))+m\sum\limits_{j=m+1}^n\Big(f(j)-n^{\alpha-1}\sum\limits_{\ell=1}^nf(\ell)\Big).$$
Thus, we have
$$\begin{array}{ll}
&(1-n^{\alpha-1}){\rm Var}_{K_n}(f)-{\rm Var}_{K_n}(M_{\alpha,K_n}f)\\
&=\displaystyle n^{\alpha-1}\Big(m(n-m-1)\sum\limits_{i=m+1}^n(f(i)-f(m+1))-\sum\limits_{m+1\leq i<j\leq n}(f(j)-f(i))\Big)\\
&\quad+\displaystyle(1-n^{\alpha-1})\sum\limits_{1\leq i<j\leq m}(f(j)-f(i))\\
&\quad+\displaystyle(n-m)(n^\alpha-1)\sum\limits_{i=1}^mf(i)+(n-m)(n-m-1)n^{\alpha-1}\sum\limits_{i=1}^m(f(m+1)-f(i))\geq 0.
\end{array}$$
So (4.2) holds.

On the other hand, it is clear that $M_{\alpha,K_n}f(i)\leq M_{\alpha,K_n}f(j)$
for any $1\leq i\leq j\leq n$. By (4.2) and Minkowski's inequality, one gets
$$\begin{array}{ll}
&{\rm Var}_{q,K_n}(M_{\alpha,K_n}f)=\displaystyle\Big(\sum\limits_{1\leq i<j\leq n}|M_{\alpha,K_n}f(j)-M_{\alpha,K_n}f(i)|^q\Big)^{1/q}\\
&\qquad\qquad\qquad\quad=\displaystyle\Big(\sum\limits_{1\leq i<j\leq n}\Big|\sum\limits_{r=i}^{j-1}(M_{\alpha,K_n}f(r+1)-M_{\alpha,K_n}f(r))\Big|^q\Big)^{1/q}\\
&\qquad\qquad\qquad\quad\leq\displaystyle\sum\limits_{r=1}^{n-1}(r(n-r))^{1/q}(M_{\alpha,K_n}f(r+1)-M_{\alpha,K_n}f(r))\\
&\qquad\qquad\qquad\quad\leq\displaystyle(n-1)^{-\alpha}\sum\limits_{r=1}^{n-1}r(n-r)(M_{\alpha,K_n}f(r+1)-M_{\alpha,K_n}f(r))\\
&\qquad\qquad\qquad\quad\leq(n-1)^{-\alpha}{\rm Var}_{K_n}(M_{\alpha,K_n}f)\leq (n-1)^{-\alpha}(1-n^{\alpha-1}){\rm Var}_{K_n}(f),
\end{array}$$
where we have used the fact $r(n-r)\ge n-1$ and $1/q-1=-\alpha<0$. This
proves (4.1) and completes the proof of Theorem \ref{thm1.5}.
\end{proof}

\begin{proof}[{\rm Proof of Theorem \ref{thm1.6}}]
Let $0<\alpha<1$ and $q=\frac{1}{1-\alpha}$. The proof will be divided
into two steps:

Without loss of generality we may assume $S_n=(V_{S_n},E_{S_n})$, where
$V_{S_n}=\{1,2,\ldots,n\}$ and $E_{S_n}=\{1\sim 2,\ldots,1\sim n\}$. It
is clear that $M_{\alpha,S_n}\delta_2(2)=1$, $M_{\alpha,S_n}\delta_2(i)
=n^{\alpha-1}$ for all $i\in\{1,3,\ldots,n\}$. It follows that
$${\bf C}_{\alpha,S_n}^{1,q}\geq\frac{{\rm Var}_{q,S_n}(M_{\alpha,S_n}f)}{{\rm Var}_{S_n}(f)}\geq((1-n^{\alpha-1})^q)^{1/q}=1-n^{\alpha-1}.$$
Thus, for (1.6) it suffices to show that for any $f:\mathbb{R}\rightarrow\mathbb{R}$
satisfying ${\rm Var}_{S_n}(f)<\infty$,
$${\rm Var}_{q,S_n}(M_{\alpha,S_n}f)\leq (1-n^{\alpha-1}){\rm Var}_{S_n}(f).\eqno(4.3)$$

Let us fix $f:\mathbb{R}\rightarrow\mathbb{R}$ satisfying
${\rm Var}_{S_n}(f)<\infty$. Without loss of generality we may assume
$f\geq0$. Since $q\geq1$, then we have ${\rm Var}_{q,S_n}(M_{\alpha,S_n}f)
\leq{\rm Var}_{S_n}(M_{\alpha,S_n}f)$. For (4.3) it is enough to show that
$${\rm Var}_{S_n}(M_{\alpha,S_n}f)\leq (1-n^{\alpha-1}){\rm Var}_{S_n}(f).\eqno(4.4)$$

Let $d_i=f(1)-f(i)$ for any $i\in\{2,\ldots,n\}$. Set $m=\sum_{i=1}^nf(i)$.
We can write
$$M_{\alpha,S_n}f(1)=\max\{f(1),n^{\alpha-1}m\},$$
$$M_{\alpha,S_n}f(i)=\max\{f(i),2^{\alpha-1}(f(1)+f(i)),n^{\alpha-1}m\},\ \ \ i\in\{2,\ldots,n\}.$$
The proof of (4.4) will be divided into two steps:

{\bf Step 1.} Proof of (4.4) for the case $f(1)>n^{\alpha-1}m$.

In this case we can write
$${\rm Var}_{S_n}(M_{\alpha,S_n}f)=\sum\limits_{i=2}^n|M_{\alpha,S_n}f(i)-f(1)|.$$
Clearly,
$$M_{\alpha,S_n}f(i)-f(1)=\max\{-d_i,(2^{\alpha}-1)f(1)-2^{\alpha-1}d_i,n^{\alpha-1}m-f(1)\},\ \ \ i\in\{2,\ldots,n\}.$$
Let $i\in\{2,\ldots,n\}$. We shall prove that
$$\begin{array}{ll}
&|\max\{-d_i,(2^{\alpha}-1)f(1)-2^{\alpha-1}d_i,n^{\alpha-1}m-f(1)\}|\\
&\leq\displaystyle(1-n^{\alpha-1})|d_i|+(2^\alpha-1)f(1)+\frac{2^\alpha-1}{n^\alpha-1}(f(1)-n^{\alpha-1}m)-\frac{(n^\alpha-n^{\alpha-1})(2^\alpha-1)}{n^\alpha-1}d_i.
\end{array}\eqno(4.5)$$
If $d_i\leq 0$, then we have
$$\max\{-d_i,(2^{\alpha}-1)f(1)-2^{\alpha-1}d_i,n^{\alpha-1}m-f(1)\}=\max\{-d_i,(2^{\alpha}-1)f(1)-2^{\alpha-1}d_i\}.$$
It is not difficult to check that $\frac{(n^\alpha-n^{\alpha-1})(2^\alpha-1)}{n^\alpha-1}
\geq n^{\alpha-1}$. This yields (4.5) in this case. If $0<d_i\leq(2-2^{1-\alpha})f(1)$,
then $(2^{\alpha}-1)f(1)-2^{\alpha-1}d_i\geq0$ and
$$\max\{-d_i,(2^{\alpha}-1)f(1)-2^{\alpha-1}d_i,n^{\alpha-1}m-f(1)\}=(2^{\alpha}-1)f(1)-2^{\alpha-1}d_i.$$
A simple computation gives that
$$1-n^{\alpha-1}-\frac{(n^\alpha-n^{\alpha-1})(2^\alpha-1)}{n^\alpha-1}\geq-2^{\alpha-1}.$$
So (4.5) holds in this case. If $(2-2^{1-\alpha})f(1)<d_i\leq f(1)$, then
$(2^{\alpha}-1)f(1)-2^{\alpha-1}d_i<0$ and
$$|\max\{-d_i,(2^{\alpha}-1)f(1)-2^{\alpha-1}d_i,n^{\alpha-1}m-f(1)\}|=\min\{2^{\alpha-1}d_i-(2^{\alpha}-1)f(1),f(1)-n^{\alpha-1}m\}.$$
We want to show that
$$\begin{array}{ll}
&\displaystyle\Big(1-n^{\alpha-1}-\frac{(n^\alpha-n^{\alpha-1})(2^\alpha-1)}{n^\alpha-1}\Big)d_i
+(2^\alpha-1)f(1)+\frac{2^\alpha-1}{n^\alpha-1}(f(1)-n^{\alpha-1}m)\\
&-\min\{2^{\alpha-1}d_i-(2^{\alpha}-1)f(1),f(1)-n^{\alpha-1}m\}\geq0.
\end{array}\eqno(4.6)$$
For convenience, we set
$$\Phi(t)=\frac{2^\alpha-1}{n^\alpha-1}t-\min\{2^{\alpha-1}d_i-(2^{\alpha}-1)f(1),t\},\ \ t\in[0,\infty).$$
Since $0<\frac{2^\alpha-1}{n^\alpha-1}<1$, it is not difficult to
see that $\Phi(t)\geq\Phi(2^{\alpha-1}d_i-(2^{\alpha}-1)f(1))$ for
all $t\geq0$. Thus, we have
$$\begin{array}{ll}
&\displaystyle\Big(1-n^{\alpha-1}-\frac{(n^\alpha-n^{\alpha-1})(2^\alpha-1)}{n^\alpha-1}\Big)d_i
+(2^\alpha-1)f(1)+\frac{2^\alpha-1}{n^\alpha-1}(f(1)-n^{\alpha-1}m)\\
&-\min\{2^{\alpha-1}d_i-(2^{\alpha}-1)f(1),f(1)-n^{\alpha-1}m\}\\
&\geq\displaystyle \Big(1-n^{\alpha-1}-\frac{(n^\alpha-n^{\alpha-1})(2^\alpha-1)}{n^\alpha-1}\Big)d_i\\
&\quad+\displaystyle(2^\alpha-1)f(1)+\Big(\frac{2^\alpha-1}{n^\alpha-1}-1\Big)(2^{\alpha-1}d_i-(2^{\alpha}-1)f(1)).
\end{array}$$
Thus, for (4.6) it suffices to show that
$$\begin{array}{ll}
&\displaystyle\Big(1-n^{\alpha-1}-\frac{(n^\alpha-n^{\alpha-1})(2^\alpha-1)}{n^\alpha-1}\Big)d_i\\
&\quad+\displaystyle(2^\alpha-1)f(1)+\Big(\frac{2^\alpha-1}{n^\alpha-1}-1\Big)(2^{\alpha-1}d_i-(2^{\alpha}-1)f(1))\geq0,
\end{array}\eqno(4.7)$$
for $d_i=(2-2^{1-\alpha})f(1)$ and $d_i=f(1)$.

When $d_i=f(1)$, we have
$$\begin{array}{ll}
&\displaystyle\Big(1-n^{\alpha-1}-\frac{(n^\alpha-n^{\alpha-1})(2^\alpha-1)}{n^\alpha-1}\Big)d_i\\
&\quad+\displaystyle(2^\alpha-1)f(1)+\Big(\frac{2^\alpha-1}{n^\alpha-1}-1\Big)(2^{\alpha-1}d_i-(2^{\alpha}-1)f(1))\\
&=\displaystyle\Big(1-n^{\alpha-1}-\frac{(n^\alpha-n^{\alpha-1})(2^\alpha-1)}{n^\alpha-1}\Big)f(1)\\
&\quad+\displaystyle(2^\alpha-1)f(1)+\Big(\frac{2^\alpha-1}{n^\alpha-1}-1\Big)(2^{\alpha-1}f(1)-(2^{\alpha}-1)f(1))\\
&=\displaystyle(2^{\alpha-1}-n^{\alpha-1})\Big(1-\frac{2^\alpha-1}{n^\alpha-1}\Big)f(1)\geq0.
\end{array}$$
When $d_i=(2-2^{1-\alpha})f(1)$, we have
$$\begin{array}{ll}
&\displaystyle\Big(1-n^{\alpha-1}-\frac{(n^\alpha-n^{\alpha-1})(2^\alpha-1)}{n^\alpha-1}\Big)d_i\\
&\quad+\displaystyle(2^\alpha-1)f(1)+\Big(\frac{2^\alpha-1}{n^\alpha-1}-1\Big)(2^{\alpha-1}d_i-(2^{\alpha}-1)f(1))\\
&=\displaystyle\Big(1+2^{\alpha-1}-n^{\alpha-1}-\frac{(n^\alpha-n^{\alpha-1})(2^\alpha-1)}{n^\alpha-1}\Big)(2-2^{1-\alpha})f(1).
\end{array}$$
Note that
$$\begin{array}{ll}
&\displaystyle 1+2^{\alpha-1}-n^{\alpha-1}-\frac{(n^\alpha-n^{\alpha-1})(2^\alpha-1)}{n^\alpha-1}\\
&=\displaystyle 1+2^{\alpha-1}-n^{\alpha-1}-(2^\alpha-1)-\frac{(2^\alpha-1)(1-n^{\alpha-1})}{n^\alpha-1}\\
&=\displaystyle(1-2^{\alpha-1})+(1-n^{\alpha-1})\Big(1-\frac{2^\alpha-1}{n^\alpha-1}\Big)\geq0.
\end{array}$$
So (4.7) holds.

Then by (4.5), we have
$$\begin{array}{ll}
&{\rm Var}_{S_n}(M_{\alpha,S_n}f)\leq\displaystyle(1-n^{\alpha-1})\sum\limits_{i=2}^m|d_i|-\frac{(n^\alpha-n^{\alpha-1})(2^\alpha-1)}{n^\alpha-1}\sum\limits_{i=2}^nd_i\\
&\qquad\qquad\qquad\quad+\displaystyle(2^\alpha-1)f(1)+\frac{2^\alpha-1}{n^\alpha-1}(f(1)-n^{\alpha-1}m)\\
&\qquad\qquad\qquad\leq(1-n^{\alpha-1}){\rm Var}_{S_n}(f),
\end{array}$$
since
$$-\frac{(n^\alpha-n^{\alpha-1})(2^\alpha-1)}{n^\alpha-1}\sum\limits_{i=2}^nd_i+(2^\alpha-1)f(1)+\frac{2^\alpha-1}{n^\alpha-1}(f(1)-n^{\alpha-1}m)=0.$$

{\bf Step 2.} Proof of (4.4) for the case $f(1)\leq n^{\alpha-1}m$.

In this case we have
$$|M_{\alpha,S_n}f(i)-M_{\alpha,S_n}f(1)|=\max\{f(i)-n^{\alpha-1}m,2^{\alpha-1}(f(1)+f(i))-n^{\alpha-1}m,0\},\ \ \ i\in\{2,\ldots,n\}.$$
We want to show that
$$\sum\limits_{i=2}^n\max\{f(i)-n^{\alpha-1}m,2^{\alpha-1}(f(1)+f(i))-n^{\alpha-1}m,0\}\leq(1-n^{\alpha-1})\sum\limits_{i=2}^n|f(i)-f(1)|.\eqno(4.7)$$
For convenience, we define the vector-valued function
$$\Gamma(\vec{u})=\sum\limits_{i=2}^n\max\{u_i-n^{\alpha-1}m(\vec{u}),2^{\alpha-1}(u_1+u_i)-n^{\alpha-1}m(\vec{u}),0\},\ \ \ u=(u_1,\ldots,u_n).$$
Here $m(\vec{u})=\sum_{i=1}^nu_i$. Let $\vec{f}=(f(1),\ldots,f(n)$.
For (4.7) it suffices to show that
$$\Gamma(\vec{f})\leq (1-n^{\alpha-1})\sum\limits_{i=2}^n|f(i)-f(1)|.\eqno(4.8)$$
It is not difficult to check that for any $\vec{u},\,\vec{v}\in\mathbb{R}^n$
and $t\geq0$,
$$\Gamma(t\vec{u})=t\Gamma(\vec{u}),\ \ \ \ \Gamma(\vec{u}+\vec{v})=2\Gamma\Big(\frac{\vec{u}+\vec{v}}{2}\Big)\leq\Gamma(\vec{u})+\Gamma(\vec{v}).\eqno(4.9)$$
For any $i\in\{1,\ldots,n\}$ we denote by $\vec{e}_i=(0,\ldots,0,1,0,
\ldots,0)\in\mathbb{R}^n$ the canonical $i$-th base vector. Applying
(4.9), one gets
$$\Gamma(\vec{e}_i)=\max\{1-n^{\alpha-1},2^{\alpha-1}-n^{\alpha-1}\}=1-n^{\alpha-1},\ \ \ i\in\{1,\ldots,n\}.\eqno(4.10)$$
Next we consider two cases:

{\it Subcases 1. $(f(1)=0)$}.  In this case we get by (4.8) and (4.9) that
$$\Gamma(\vec{f})=\Gamma\Big(\sum\limits_{i=2}^nf(i)\vec{e}_i\Big)\leq\sum\limits_{i=2}^nf(i)\Gamma(\vec{e}_i)\leq (1-n^{\alpha-1})\sum\limits_{i=2}^n|f(i)-f(1)|,$$
which gives (4.8) in this case.

{\it Subcases 2. $(f(1)>0)$}. Let $\vec{u}=(u_1,\ldots,u_n)$ with $u_1=1$
and $u_i=0$ or $u_i=1$. Assume $|\{u_i=1:1\leq i\leq n\}|=k$. We shall
prove
$$(1-n^{\alpha-1})(n-k)-\Gamma(\vec{u})\geq0.\eqno(4.11)$$
Observe that
$$\max\{u_i-n^{\alpha-1}m(\vec{u}),2^{\alpha-1}(u_1+u_i)-n^{\alpha-1}m(\vec{u}),0\}
=\max\{2^{\alpha}-n^{\alpha-1}k,0\},\ \ {\rm if}\ \ u_i=1,$$
$$\max\{u_i-n^{\alpha-1}m(\vec{u}),2^{\alpha-1}(u_1+u_i)-n^{\alpha-1}m(\vec{u}),0\}
=\max\{2^{\alpha-1}-n^{\alpha-1}k,0\},\ \ {\rm if}\ \ u_i=0.$$
It follows that
$$\Gamma(\vec{u})\leq(k-1)\max\{2^{\alpha}-n^{\alpha-1}k,0\}+(n-k)\max\{2^{\alpha-1}-n^{\alpha-1}k,0\}.$$
If $n^{\alpha-1}k\geq 2^\alpha$, then
$$\Gamma(\vec{u})\leq(k-1)\max\{2^{\alpha}-n^{\alpha-1}k,0\}+(n-k)\max\{2^{\alpha-1}-n^{\alpha-1}k,0\}=0.$$
In this case (4.11) is clear. If $2^{\alpha-1}-n^{\alpha-1}k>0$, then we have
$$\begin{array}{ll}
&(1-n^{\alpha-1})(n-k)-\Gamma(\vec{u})=(1-n^{\alpha-1})(n-k)
-(k-1)(2^{\alpha}-n^{\alpha-1}k)-(n-k)(2^{\alpha-1}-n^{\alpha-1}k)\\
&\qquad\qquad\qquad\qquad\qquad=(n-k)(1-2^{\alpha-1})+(k-1)(n^\alpha-2^\alpha)\geq0.
\end{array}$$
If $2^{\alpha-1}\leq n^{\alpha-1}k\leq 2^\alpha$, then we have $k\geq2$ and
$$\begin{array}{ll}
&(1-n^{\alpha-1})(n-k)-\Gamma(\vec{u})=(1-n^{\alpha-1})(n-k)-(k-1)(2^{\alpha}-n^{\alpha-1}k)\\
&\qquad\qquad\qquad\qquad\qquad=n^{\alpha-1}k^2-(1+2^\alpha)k+n-n^\alpha+2^\alpha=:g(n,x).
\end{array}$$
Set $x_n=\frac{(1+2^\alpha)n^{1-\alpha}}{2}$. Since $n^{\alpha-1}>0$
and $\frac{\partial g(n,x)}{dx}=2n^{\alpha-1}x-(1+2^\alpha)$, then we have
$$g(n,x)\geq g(n,2)=4n^{\alpha-1}-2(1+2^\alpha)+n-n^\alpha+2^\alpha=:h_1(n),\ \ {\rm if}\ x_n\leq2;$$
$$g(n,x)\geq g(n,x_n)=n-n^\alpha+2^\alpha-\frac{1}{4}(1+2^\alpha)^2n^{1-\alpha}=:h_2(n),\ \ {\rm if}\ x_n>2.$$
Since $h_1'(n)=4(\alpha-1)n^{\alpha-2}+1-\alpha n^{\alpha-1}\geq0$ and
$h_2'(n)=1-\alpha n^{\alpha-1}-\frac{1-\alpha}{4n^\alpha}(1+2^\alpha)^2\geq0$,
then we have
$$h_1(n)\geq h_1(4)=2-2^\alpha>0,$$
$$h_2(n)\geq h_2(4)=4-4^\alpha+2^\alpha-\frac{1}{4}(1+2^\alpha)^24^{1-\alpha}=\frac{(2^\alpha-1)(1+3\cdot 2^\alpha-2^{3\alpha})}{4^\alpha}>0.$$
Thus, (4.11) holds in this case.

Now let $\vec{b}=(1,b_2,\ldots,b_n)$ with $b_i\in[0,1]$ for $i\in\{2,\ldots,n\}$.
We want to show that
$$\Gamma(\vec{b})\leq(1-n^{\alpha-1})\sum\limits_{i=2}^n(1-b_i).\eqno(4.12)$$
We can write
$$\vec{b}=\sum\limits_{\epsilon}\lambda_\epsilon\vec{w}_\epsilon,$$
where
$$\epsilon=(\epsilon_2,\ldots,\epsilon_n)\in(\{0,1\})^{n-1},\ \ \vec{w}_\epsilon=(1,\epsilon_2,\ldots,\epsilon_n),
\ \ \lambda_\epsilon=\prod\limits_{j=2}^n(\epsilon_ib_i+(1-\epsilon_i)(1-b_i),\ \ \sum\limits_{\epsilon}\lambda_\epsilon=1.$$
Applying (4.11), one gets
$$\Gamma(\vec{b})\leq\displaystyle\sum\limits_{\epsilon}\lambda_\epsilon\Gamma(\vec{w}_\epsilon)
\leq\displaystyle(1-n^{\alpha-1})\sum\limits_{\epsilon}\lambda_\epsilon\sum\limits_{i=2}^n(1-\epsilon_i)
\leq (1-n^{\alpha-1})\sum\limits_{i=2}^n(1-b_i).$$
This yields (4.12).

Finally, we write
$$\vec{f}=f(1)\vec{b}+\sum\limits_{i=2}^nr_i\vec{e}_i,$$
where
$$\vec{b}=(1,b_2,\ldots,b_n),\ \  b_i=\min\Big\{\frac{f(i)}{f(1)},1\Big\},\ \ \ \ i\in\{2,\ldots,n\},$$
$$r_i=\max\{f(i)-f(1),0\},\ \ \ \ \ i\in\{2,\ldots,n\}.$$
It follows from (4.9), (4.10) and (4.12) that
$$\Gamma(\vec{f})\leq f(1)\Gamma(\vec{b})+\sum\limits_{i=2}^nr_i\Gamma(\vec{e}_i)\leq (1-n^{\alpha-1})
\sum\limits_{i=2}^n(f(1)(1-b_i)+r_i)\leq (1-n^{\alpha-1})\sum\limits_{i=2}^n|f(i)-f(1)|,$$
which gives (4.8) in this case.
\end{proof}
\begin{proof}[{\rm Proof of Theorem \ref{thm1.7}}]
Let $n\geq4$, $\alpha\in(0,1)$ and $q=\frac{1}{1-\alpha}$. Without
loss of generality we may assume $S_n=(V_{S_n},E_{S_n})$, where
$V_{S_n}=\{1,2,\ldots,n\}$ and $E_{S_n}=\{1\sim 2,\ldots,1\sim n\}$.
Observe that
$$\widetilde{M}_{\alpha,S_n}\delta_2(1)=2^{\alpha-1},\ \ \ \widetilde{M}_{\alpha,S_n}\delta_2(2)=1,\ \ \widetilde{M}_{\alpha,S_n}\delta_2(i)=n^{\alpha-1},\ \ i\in\{3,\ldots,n\}.$$
It follows that
$$\widetilde{{\bf C}}_{\alpha,S_n}^{1,q}\geq\frac{{\rm Var}_{q,S_n}(\widetilde{M}_{\alpha,S_n}\delta_2)}{{\rm Var}_{S_n}(\delta_2)}
\geq((1-2^{\alpha-1})^q+(n-2)(2^{\alpha-1}-n^{\alpha-1})^q)^{1/q}.$$
Thus, for (1.7) it suffices to show that for any $f:V_{S_n}\rightarrow\mathbb{R}$
satisfying ${\rm Var}_{S_n}(f)>0$ such that
$${\rm Var}_{q,S_n}(\widetilde{M}_{\alpha,S_n}f)\leq((1-2^{\alpha-1})^q+(n-2)(2^{\alpha-1}-n^{\alpha-1})^q)^{1/q}{\rm Var}_{S_n}(f).\eqno(4.13)$$

In what follows, let $f:V_{S_n}\rightarrow\mathbb{R}$ satisfy
${\rm Var}_{S_n}(f)>0$. We may assume without loss of generality
$f\geq0$. For convenience, we set
$$d_i=|f(i)-f(1)|,\ \ \ i\in\{2,\ldots,n\}.$$
Then we have ${\rm Var}_{S_n}(f)=\sum_{i=2}^nd_i$ and
$$\widetilde{M}_{\alpha,S_n}f(1)=\max\Big\{f(1),2^{\alpha-1}(f(1)+\max_{2\leq j\leq n}f(j)),n^{\alpha-1}\sum_{i=1}^{n}f(i)\Big\},$$
$$\widetilde{M}_{\alpha,S_n}f(i)=\max\Big\{f(i),2^{\alpha-1}(f(1)+f(i)),n^{\alpha-1}\sum_{i=1}^{n}f(i)\Big\},\ \ \ \ i\in\{2,\ldots,n\}.$$
The proof of (4.13) will be divided into three steps:

{\bf Step 1.} At first we prove the following estimate
$$|\widetilde{M}_{\alpha,S_n}f(i)-\widetilde{M}_{\alpha,S_n}f(1)|\leq(1-2^{\alpha-1})d_i+
(2^{\alpha-1}-n^{\alpha-1})\sum\limits_{j=2,\,j\neq i}^{n}d_j,\ \ \ \ i\in\{2,\ldots,n\}.\eqno(4.14)$$
In order to prove (4.14), let us firstly prove the following estimates:
$$\min\{n^{\alpha-1},1-2^{\alpha-1}\}+(n-3)\min\{n^{\alpha-1},2^{\alpha-1}-n^{\alpha-1}\}\geq2(2^{\alpha-1}-n^{\alpha-1}),\eqno(4.15)$$
$$\min\{n^{\alpha-1},t\}f(1)\leq n^{\alpha-1}x+t|x-f(1)|,\ \ \ x,\,t\in[0,\infty),\eqno(4.16)$$
$$2^{\alpha-1}(f(1)+f(j))-n^{\alpha-1}\sum\limits_{k=1}^nf(k)\leq(1-2^{\alpha-1})d_k+(2^{\alpha-1}-n^{\alpha-1})\sum\limits_{k=2,\,k\neq i}d_k,\ \ \forall i\neq j.\eqno(4.17)$$

{\bf Proof of (4.15)}. For (4.15) it suffices to show that
$$\begin{array}{ll}
&\min\{(n-2)n^{\alpha-1},n^{\alpha-1}+(n-3)(2^{\alpha-1}-n^{\alpha-1}),1-2^{\alpha-1}+(n-3)n^{\alpha-1},\\
&\quad1-2^{\alpha-1}+(n-3)(2^{\alpha-1}-n^{\alpha-1})\}\geq2(2^{\alpha-1}-n^{\alpha-1}).
\end{array}\eqno(4.18)$$
Clearly,
$$(n-2)n^{\alpha-1}-2(2^{\alpha-1}-n^{\alpha-1})=n^\alpha-2^\alpha>0.$$
$$4^{\alpha-1}+(2^{\alpha-1}-4^{\alpha-1})-2(2^{\alpha-1}-4^{\alpha-1})=2^{\alpha-1}(2^\alpha-1)>0,$$
$$n^{\alpha-1}+(n-3)(2^{\alpha-1}-n^{\alpha-1})-2(2^{\alpha-1}-n^{\alpha-1})\geq
n^{\alpha-1}+(n-5)(2^{\alpha-1}-n^{\alpha-1})>0,\ \ n\geq5,$$
$$\begin{array}{ll}
&1-2^{\alpha-1}+(n-3)n^{\alpha-1}-2(2^{\alpha-1}-n^{\alpha-1})\\
&=(n^\alpha-n^{\alpha-1})-3\cdot 2^{\alpha-1}+1\geq 4^\alpha-4^{\alpha-1}-3\cdot 2^{\alpha-1}+1
=3(2^{\alpha-1}-2^{-1})^2+\frac{1}{4}>0,
\end{array}$$
$$1-2^{\alpha-1}-(2^{\alpha-1}-4^{\alpha-1})=(1-2^{\alpha-1})^2>0,$$
$$1-2^{\alpha-1}+(n-5)(2^{\alpha-1}-n^{\alpha-1})>0,\ \ n\geq5.$$
So (4.18) holds.

{\bf Proof of (4.16)}. If $n^{\alpha-1}\geq t$, then
$$tf(1)\leq t(x+|x-f(1)|)\leq n^{\alpha-1}x+t|x-f(1)|.$$
If $n^{\alpha-1}\leq t$, then
$$n^{\alpha-1}f(1)\leq n^{\alpha-1}(x+|x-f(1)|)\leq n^{\alpha-1}x+t|x-f(1)|.$$
So (4.16) holds.

{\bf Proof of (4.17)}. Observe that for any $i\neq j$,
$$f(1)+f(j)-d_j=2\min\{f(1),f(j)\}\leq2f(1).$$
Applying (4.16), we have
$$n^{\alpha-1}f(i)+(1-2^{\alpha-1})d_i\geq\min\{n^{\alpha-1},t\}f(1),$$
$$n^{\alpha-1}f(k)+(2^{\alpha-1}-n^{\alpha-1})d_k\geq\min\{n^{\alpha-1},2^{\alpha-1}-n^{\alpha-1}\}f(1).$$
Then we get from (4.15) that
$$\begin{array}{ll}
&\displaystyle 2^{\alpha-1}(f(1)+f(j))-n^{\alpha-1}\sum\limits_{i=1}^nf(i)-\Big((1-2^{\alpha-1})d_i+(2^{\alpha-1}-n^{\alpha-1})\sum\limits_{k=2,\,k\neq i}d_k\Big)\\
&=(2^{\alpha-1}-n^{\alpha-1})(f(1)+f(j)-d_j)-\Big(n^{\alpha-1}f(i)+(1-2^{\alpha-1})d_i\Big)\\
&\quad\displaystyle-\sum\limits_{k\in\{2,\ldots,n\}\setminus\{i,j\}}^{n}
\Big(n^{\alpha-1}f(k)+(2^{\alpha-1}-n^{\alpha-1})d_k\Big)\\
&\leq\Big(2(2^{\alpha-1}-n^{\alpha-1})-\min\{n^{\alpha-1},1-2^{\alpha-1}\}-(n-3)\min\{n^{\alpha-1},2^{\alpha-1}-n^{\alpha-1}\}\Big)f(1)\leq0,
\end{array}$$
which leads to (4.17).
	
We now turn to the proof of (4.14). Let us fix $i\in\{2,\ldots,n\}$.
If $\widetilde{M}_{\alpha,S_n}f(i)>\widetilde{M}_{\alpha,S_n}f(1)$,
then we have $\widetilde{M}_{\alpha,S_n}f(i)=f(i)>f(1)$ and
$$|\widetilde{M}_{\alpha,S_n}f(i)-\widetilde{M}_{\alpha,S_n}f(1)|
\leq \widetilde{M}_{\alpha,S_n}f(i)-\widetilde{M}_{\alpha,S_n}f(1)
\leq f(i)-2^{\alpha-1}(f(1)+f(i))\leq(1-2^{\alpha-1})d_i.$$
This proves (4.14) in this case.

If $\widetilde{M}_{\alpha,S_n}f(1)>\widetilde{M}_{\alpha,S_n}f(i)$,
then we consider two cases:

{\it Subcase 1. ($\widetilde{M}_{\alpha,S_n}f(1)=f(1)$)}. In this
case we have $f(1)>f(i)$, and
$$\widetilde{M}_{\alpha,S_n}f(1)-\widetilde{M}_{\alpha,S_n}f(i)\leq f(1)-2^{\alpha-1}(f(1)+f(i))\leq(1-2^{\alpha-1})d_i.$$

{\it Subcase 2. ($\widetilde{M}_{\alpha,S_n}f(1)\neq f(1)$)}. In this
case there exists $j_0\in\{2,\ldots,n\}$ such that $j_0\neq i$,
$f(j_0)=\max_{2\leq j\leq n}f(j)$ and
$$\widetilde{M}_{\alpha,S_n}f(1)=2^{\alpha-1}(f(1)+f(j_0)).$$
Then we get from (4.17) that
$$\begin{array}{ll}
&\widetilde{M}_{\alpha,S_n}f(1)-\widetilde{M}_{\alpha,S_n}f(i)\leq
2^{\alpha-1}(f(1)+f(j_0))-n^{\alpha-1}\sum\limits_{j=1}^{n}f(j)\\
&\qquad\qquad\qquad\qquad\quad\leq\displaystyle(1-2^{\alpha-1})d_i+(2^{\alpha-1}-n^{\alpha-1})\sum\limits_{j=2,\,j\neq i}^{n}d_i.
\end{array}$$
This proves (4.14).

{\bf Step 2. Proof of (4.13)}. By (4.14) and Minkowski's inequality, one gets

$$\begin{array}{ll}
&{\rm Var}_{q,S_n}(\widetilde{M}_{\alpha,S_n}f)\leq\displaystyle\Big(\sum\limits_{i=2}^{n}\Big((1-2^{\alpha-1})d_i+(2^{\alpha-1}-n^{\alpha-1})
			\sum\limits_{j=2,\,j\neq i}^{n}d_j\Big)^q\Big)^{1/q}\\
&\qquad\qquad\qquad\quad\displaystyle\leq\sum\limits_{j=2}^{n}\Big(((1-2^{\alpha-1})d_j)^q+
			\sum\limits_{i=2,\,i\neq j}^{n}((2^{\alpha-1}-n^{\alpha-1})d_j)^q\Big)^{1/q}\\
&\qquad\qquad\qquad\quad\leq\displaystyle\Big((1-2^{\alpha-1})^q+
			(n-2)(2^{\alpha-1}-n^{\alpha-1})^q\Big)^{1/q}\sum\limits_{j=2}^{n}d_j\\
&\qquad\qquad\qquad\quad\leq\displaystyle\Big((1-2^{\alpha-1})^q+
			(n-2)(2^{\alpha-1}-n^{\alpha-1})^q\Big)^{1/q}{\rm Var}_{S_n}(f).
\end{array}$$
Here we used $q>1$ and $d_j\geq0$. This proves (4.13) and finishes the proof
of Theorem \ref{thm1.7}.
	
\end{proof}

\bigskip

\quad\hspace{-20pt}{\bf Funding}

\quad\hspace{-20pt}Feng Liu is supported by the Natural Science Foundation
of Shandong Province (Grant No. ZR2023MA022).

\medskip

\quad\hspace{-20pt}{\bf Competing Interests}

\quad\hspace{-20pt}The authors declare that there is no competing interest
regrading the publication of this paper.


\begin{thebibliography}{99}

\bibitem{Ald}{J. M. Aldaz},
{The weak type $(1,\,1)$ bounds for the maximal function associated to cubes grow
to infinity with the dimension}, Ann. of Math. {\bf 173} (2011), 1013--1023.


\bibitem{AL}{J. M. Aldaz and J. P\'{e}rez L\'{a}zaro},
{Functions of bounded variation, the derivative of the one dimensional maximal
function, and applications to inequalities}, Trans. Am. Math. Soc. {\bf 359} (2007),
2443--2461.

\bibitem{BM}{N. Badr and J. M. Martell}, {Weighted norm inequalities on graphs},
J. Geom. Anal. {\bf 22} (2012), 1173--1210.

\bibitem{BCHP}J. Bober, E. Carneiro, K. Hughes and L. B. Pierce,
On a discrete version of Tanaka's theorem for maximal functions,
Proc. Am. Math. Soc. {\bf 140} (2012), 1669--1680.

\bibitem{CFS}E. Carneiro, R. Finder and M. Sousa,
{On the variation of maximal operators of convolution type II},
Rev. Mat. Iberoam. {\bf 34} (2018), 739--766.

\bibitem{CM}E. Carneiro and J. Madrid,
{Derivative bounds for fractional maximal functions},
Trans. Am. Math. Soc. {\bf 369} (2017), 4063--4092.


\bibitem{CS}E. Carneiro and B. F. Svaiter,
On the variation of maximal operators of convolution type,
J. Funct. Anal. {\bf 265} (2013), 837--865.

\bibitem{CMS}M. Cowling, S. Meda and A. G. Setti,
{Estimates for functions of the Laplace operator on homogeneous trees},
{Trans. Amer. Math. Soc.} {\bf 352}(9) (2000), 4271--4293.

\bibitem{GKM}C. Gonz\'{a}lez-Riquelme, V. Kovac and J. Madrid,
Sharp variational inequalities for the Hardy--Littlewood maximal operator
on finite undirected graphs, arXiv: 2603.12462v1.


\bibitem{GM}C. Gonz\'{a}lez-Riquelme and J. Madrid,
Sharp inequalities for maximal operators on finite graphs,
J. Geom. Anal. {\bf 31}(10) (2021), 9708--9744.

\bibitem{GM2}C. Gonz\'{a}lez-Riquelme and J. Madrid,
Sharp inequalities for maximal operators on finite graphs, II,
J. Math. Anal. Appl. {\bf 506} (2022), 125647.

\bibitem{GM3}{L. Grafakos and S. Montgomery-Smith},
{Best constants for uncentred maximal functions},
Bull. London Math. Soc. {\bf 29}(1) (1997), 60--64.

\bibitem{GMM}{L. Grafakos, S. Montgomery-Smith and O. Motrunich},
{A sharp estimate for the Hardy-Littlewood maximal function},
Studia Math. {\bf 134}(1) (1999), 57--67.

\bibitem{KP}A. Kor\'{a}nyi and M. A. Picardello,
{Boundary behaviour of eigenfunctions of the Laplace operator on trees},
{Ann. Scuola Norm. Sup. Pisa Cl. Sci. (4)} {\bf 13}(3) (1986), 389--399.

\bibitem{Ku}O. Kurka, On the variation of the Hardy--Littlewood maximal
function, Ann. Acad. Sci. Fenn. Math. {\bf 40} (2015), 109--133.

\bibitem{LX}F. Liu and Q. Xue, On the variation of the Hardy--Littlewood maximal functions on finite graphs,
Collect. Math. {\bf 72}(2) (2021), 333--349.

\bibitem{Me1}{A. D. Melas},
{On the centered Hardy-Littlewood maximal operator},
Trans. Amer. Math. Soc. {\bf 354}(8) (2002), 3263--3273.

\bibitem{Me2}{A. D. Melas},
{The best constant for the centered Hardy-Littlewood maximal inequality},
Ann. of Math. (2) {\bf 157}(2) (2003), 647--688.

\bibitem{Ra}{J. P. Ramos}, Sharp total variation results for maximal functions,
Ann. Acad. Sci. Fenn. Math. {\bf 44} (2019), 41--64.

\bibitem{ST}J. Soria and P. Tradacete,
{Best constants for the Hardy-Littlewood maximal operator on finite graphs},
{J. Math. Anal. Appl.} {\bf 436}(2) (2016), 661--682.

\bibitem{ST2}J. Soria and P. Tradacete,
{Geometric properties of infinte graphs and the Hardy-Littlewood maximal operator},
{J. Anal. Math.} {\bf 137}(2) (2019), 913--937.
		
\bibitem{XZ}C. Xiang and F. Zhao,
Sharp constants for fractional Hardy-Littlewood maximal operators on finite graphs,
Forum Math. {\bf 38}(2) (2026), 411--434.


\bibitem{Wei1}J. Weigt,
Variation of the uncentered maximal characteristic function,
Rev. Mat. Iberoam. {\bf 38} (2022), 823--849.

\bibitem{Wei2}J. Weigt, {Variation of the dyadic maximal function},
Int. Math. Res. Not. {\bf 2023} (2023), 3576--3596.


\bibitem{Wei3}J. Weigt,
{The variation of the uncentered maximal operator with respect to cubes},
J. European Math. Soc. {\bf 27}(6) (2024), 2531--2570.







\end{thebibliography}
\end{document}